\documentclass[12pt]{amsart}

\usepackage{amsthm}
\usepackage{amssymb}
\usepackage{amsmath}
\usepackage{amsfonts} 
\usepackage{times}

\newcommand{\excise}[1]{}
\newcommand{\comment}[1]{{$\star$\sf\textbf{#1}$\star$}}

\newtheorem{thm}{Theorem}[section]
\newtheorem{lemma}[thm]{Lemma}

\newtheorem{cor}[thm]{Corollary}
\newtheorem{prop}[thm]{Proposition}

\newtheorem{porism}[thm]{Porism}

\theoremstyle{definition}
\newtheorem{defn}[thm]{Definition}
\newtheorem{example}[thm]{Example}
\newtheorem{remark}[thm]{Remark}

\newtheorem{notn}[thm]{Notation}

\newenvironment{numbered}%
        {\begin{list}
                {\noindent\makebox[0mm][r]{\arabic{enumi}.}}
                {\leftmargin=5.5ex \usecounter{enumi}}
        }
        {\end{list}}

        {\begin{list}
                {\noindent\makebox[0mm][r]{(\roman{enumi})}}
                {\leftmargin=5.5ex \usecounter{enumi}}
        }
        {\end{list}}

 \newcounter{separated}
 \newcounter{starredcounter}
\newenvironment{stareqn}%
	{
	 \setcounter{separated}{\value{equation}}
	 \setcounter{equation}{\value{starredcounter}}
	 \begin{eqnarray}
	}
	{\end{eqnarray}%
	 \stepcounter{starredcounter}%
	 \setcounter{equation}{\value{separated}}%
	}

\newcommand{\baseRing}[1]{\ensuremath{\mathbb{#1}}}
\newcommand{\CC}{\baseRing{C}}
\newcommand{\NN}{\baseRing{N}}
\newcommand{\RR}{\baseRing{R}}
\newcommand{\ZZ}{\baseRing{Z}}

\def\<{\langle}
\def\>{\rangle}

\def\DD{{\mathbb D}}

\def\calE{{\mathcal E}}
\def\calF{{\mathcal F}}
\def\calG{{\mathcal G}}
\def\calH{{\mathcal H}}
\def\calK{{\mathcal K}}
\def\calM{{\mathcal M}}

\def\calO{{\mathcal O}}

\def\frakm{{\mathfrak m}}
\def\frakp{{\mathfrak p}}

\def\pp{{\mathfrak p}}

\def\boldy{{\mathbf y}}

\def\th{{\rm th}}

\def\dlim#1{{\lim\limits_{\aoverb\longrightarrow{#1}}}}
\def\ext{\operatorname{Ext}}
\def\hom{\operatorname{Hom}}
\def\gr{\operatorname{gr}}
\def\qdeg{\operatorname{qdeg}}
\def\rk{\operatorname{rank}}
\def\sol{\operatorname{\mathcal S\hspace{-.23ex}\mathit{ol}}}
\def\spec{\operatorname{Spec}}

\def\Tor{\tor}
\def\tor{\operatorname{Tor}}
\def\tot{\operatorname{Tot}}
\def\tdeg{\operatorname{tdeg}}

\def\vol{\operatorname{vol}}

\def\blank{\hspace{.1ex}{\raisebox{.2ex}{\underline{\ \,}}\hspace{.18ex}}}
\def\0{\mathbf{0}}
\def\AF{{A_F}}
\def\IAF{I_A^F}
\def\IN{\operatorname{\mathsf{in}}}
\def\SA{S_{\hspace{-.2ex}A}}
\def\SF{S_{\hspace{-.2ex}F}}
\def\SFp{S_{\hspace{-.2ex}F'}}
\def\SFk{S_{\hspace{-.2ex}F_k}}
\def\del{\partial}
\def\iff{\Leftrightarrow}

\def\vea{\varepsilon_{\hspace{-.1ex}A}}
\def\MAFB{\calM_{\raisebox{-.06ex}{$\scriptstyle\beta$}}%
	{}^{\hspace{-1.1ex}\raisebox{.06ex}{$\scriptstyle\AF$}}}
\def\MAFBp{\calM_{\raisebox{-.06ex}{$\scriptstyle\beta$}}%
	{}^{\hspace{-1.1ex}\raisebox{.06ex}{$\scriptstyle A_F'$}}}

\def\from{\leftarrow}
\def\into{\hookrightarrow}
\def\onto{\twoheadrightarrow}
\def\spot{{\hbox{\raisebox{1.7pt}{\large\bf .}}\hspace{-.5pt}}}
\def\minus{\smallsetminus}

\def\ol#1{{\overline {#1}}}

\def\ideal#1{\langle#1\rangle}
\newcommand{\aoverb}[2]{{\genfrac{}{}{0pt}{1}{#1}{#2}}}

\begin{document}

\title{Homological methods for hypergeometric families}

\author{Laura Felicia Matusevich}
\address{Department of Mathematics, Harvard University, Cambridge MA}
\email{laura@math.harvard.edu}
\thanks{LFM was partially supported by postdoctoral fellowships from
	MSRI and the NSF}

\author{Ezra Miller}
\address{School of Mathematics, University of Minnesota, Minneapolis MN}
\email{ezra@math.umn.edu}
\thanks{EM was partially supported by the NSF}

\author{Uli Walther}
\address{Department of Mathematics, Purdue University, West Lafayette IN}
\email{walther@math.purdue.edu}
\thanks{UW was partially supported by the NSF, the DfG and the
	Humboldt foundation}

\date{17 June 2004}\mbox{}

\begin{abstract}
We analyze the behavior of the holonomic rank in families of holonomic
systems over complex algebraic varieties by providing homological
criteria for rank-jumps in this general setting.  Then we investigate
rank-jump behavior for hypergeometric systems~$H_A(\beta)$ arising
from a $d \times n$ integer matrix~$A$ and a parameter $\beta \in
\CC^d$.  
To do so we introduce an Euler--Koszul functor for hypergeometric families
over~$\CC^d$, whose homology generalizes the notion of a hypergeometric
system, and we prove a homology isomorphism with our general homological
construction above.
We show that a parameter $\beta \in \CC^d$ is rank-jumping
for $H_A(\beta)$ if and only if $\beta$ lies in the Zariski closure of
the set of $\ZZ^d$-graded degrees~$\alpha$ where the local cohomology
$\bigoplus_{i < d} H^i_\frakm(\CC[\NN A])_\alpha$ of the semigroup
ring $\CC[\NN A]$ supported at its maximal graded ideal~$\frakm$ is
nonzero.  Consequently, $H_A(\beta)$ has no rank-jumps over~$\CC^d$ if
and only if $\CC[\NN A]$ \mbox{is Cohen--Macaulay of dimension~$d$}.
\end{abstract}

\maketitle
\vspace{-1ex}
\setcounter{tocdepth}{1}
\tableofcontents

\section{Introduction}

\label{sec:intro}

In the late 1980s, Gelfand, Graev and Zelevinsky introduced a class of
systems of linear partial differential equations closely related to
toric varieties \cite{GGZ87}.  These systems, now called {\em
GKZ\/~systems}, or {\em $A$-hypergeometric systems}, are constructed
from a $d\times n$ integer matrix~$A$ of rank~$d$ and a complex
parameter vector $\beta\in \CC^d$, and are denoted by~$H_A(\beta)$.
Although \mbox{$A$-hypergeometric} systems were originally introduced
with a view toward representation theory, they are well suited to
algebraic geometry in various applications.  For example, solutions of
\mbox{$A$-hypergeometric} systems appear as toric residues and as the
generating functions for intersection numbers on moduli spaces of
curves.  In another instance, the Picard--Fuchs equations governing
the variation of Hodge structures for Calabi--Yau toric hypersurfaces
are very special cases of $A$-hypergeometric systems.

The first fundamental results about the systems~$H_A(\beta)$ were
proved by Gelfand, Graev, Kapranov, and Zelevinsky.  These results
concerned the case where the semigroup~$\NN A$ generated by the
columns of~$A$ gives rise to a semigroup ring~$\CC[\NN A]$ that is
Cohen--Macaulay and graded in the standard $\ZZ$-grading \cite{GGZ87,
GKZ89}.  In geometric terms, the associated toric variety~$X_A$ is
projective and arithmetically Cohen--Macaulay.  The above authors
showed that in this case, the system~$H_A(\beta)$ gives a holonomic
module over the ring~$D$ of polynomial $\CC$-linear differential
operators in $n$ variables, and hence $H_A(\beta)$ has finite {\em
rank}\/; that is, the dimension of its space of holomorphic solutions
is finite.  Furthermore, they showed that this dimension can be
expressed combinatorially, as the integer~$\vol(A)$ that is $d!$ times
the Euclidean volume of the convex hull of the columns of $A$ and the
origin $0 \in \ZZ^d$.
The remarkable fact is that their rank formula held independently of
the parameter~$\beta$.

Even if $\CC[\NN A]$ is not Cohen--Macaulay or $\ZZ$-graded, Adolphson
showed in \cite{Adolphson-Duke94} that $H_A(\beta)$ is always a
holonomic ideal.  He further proved that, for all parameters outside
of a closed locally finite arrangement of countably many
`semi-resonant' affine hyperplanes, the characterization of rank
through volume is still correct.

It came as quite a surprise when in \cite{Sturmfels-Takayama} an
example was given showing that if $\CC[\NN A]$ is not Cohen--Macaulay
then not all parameters $\beta$ have to give the same rank.  One is
hence prompted to introduce $\calE_A$, the collection of those {\em
exceptional parameters} $\beta\in\CC^d$ for which the rank does not
take the expected value.  Nearly at the same time the case of
projective toric curves was discussed completely in \cite{CAD}: the
set~$\calE_A$ of exceptional parameters is finite in this case, and
empty precisely when $\CC[\NN A]$ is Cohen--Macaulay; moreover, at
each $\beta\in\calE_A$ the rank exceeds the volume by exactly~$1$.  It
was shown soon after in \cite{SST} that the rank can never be smaller
than the volume as long as $\CC[\NN A]$ is $\ZZ$-graded, and it was
established in the same book that $\calE_A$ is in fact contained in a
finite affine subspace arrangement.  More recently the much stronger
fact emerged that $\calE_A$ is a finite union of Zariski locally
closed sets by means of Gr\"obner basis techniques \cite{Laura-exc}.
While rank-jumps can be arbitrarily large \cite{MW}, the absence of
rank-jumping parameters is equivalent to the Cohen--Macaulay property
for $\ZZ$-graded $\CC[\NN A]$ when either $\CC[\NN A]$ has codimension
two \cite{Laura-codim2}, or if the convex hull of $A$ is a simplex
\cite{Saito-Duke}, or if $\CC[\NN A]$ is a polynomial ring modulo a
generic toric ideal \cite{Laura-exc}.

Encouraged by these results, which suggest an algebraic structure on
the set of exceptional parameters, it was conjectured in \cite{MM}
that the obstructions to the Cohen--Macaulayness of $\CC[\NN A]$ and
the set of exceptional parameters are identified in an explicit
manner.  To be precise, let $H^{< d}_\frakm(\CC[\NN A])$ be the direct
sum of all the local cohomology modules supported at the maximal
homogeneous ideal~$\frakm$ of $\CC[\NN A]$ in cohomological degrees
less than~$d$.  Then define the set~$E_A$ of {\em exceptional
quasi-degrees}\/ to be the Zariski closure in~$\CC^d$ of the set of
$\ZZ^d$-graded degrees~$\alpha$ such that $H^{< d}_\frakm(\CC[\NN A])$
has a nonzero element in degree~$-\alpha$.  With this notation, the
motivating result in this article is the following.

\begin{thm} \label{conj:MM}
For any rank $d$ matrix $A\in\ZZ^{d\times n}$ the set~$\calE_A$ of
exceptional (that is, rank-jumping) parameters is identical to the
set~$E_A$ of exceptional quasi-degrees.
\end{thm}

We note that there is no assumption on $\CC[\NN A]$ being $\ZZ$-graded.
The $\ZZ$-graded simplicial case of this result was proved in
\cite{MM} using results of \cite{Saito-Duke}.


\subsection*{Methods and results}
The first step in our proof of Theorem~\ref{conj:MM} is to construct a
homological theory to systematically detect rank-jumps.  In
Sections~\ref{sec:upper} and~\ref{sec:criteria}  we study rank
variation in any family of holonomic modules over any base~$B$, and
not just $A$-hypergeometric families over $B = \CC^d$.  The idea is
that under a suitable coherence assumption (Definition~\ref{d:hol}),
holonomic ranks behave like fiber dimensions in families of algebraic
varieties.  In particular, rank is constant almost everywhere and can
only increase on closed subsets of~$B$ (Theorem~\ref{t:semicont}).  We
develop a computational tool to check for rank-jumps at a smooth point
$\beta \in B$: since the rank-jump occurs through a failure of
flatness at~$\beta$, ordinary Koszul homology detects it
(Theorem~\ref{t:cm} and Corollary~\ref{c:cm}).

The second step toward Theorem~\ref{conj:MM} is to construct a
homolo\-gical theory for $D$-modules that reproduces the set~$E_A$ of
exceptional quasi-degrees, which a~priori arises from the commutative
notion of local cohomology.  Our main observation along these lines is
that the {\em Euler--Koszul complex}, which was already known to
Gelfand, Kapranov, and Zelevinsky for Cohen--Macaulay $\ZZ$-graded
semigroup rings \cite{GKZ89}, generalizes to fill this need.
Adolphson \cite{Adolphson-Rend} recognized that when the semigroup is
not Cohen--Macaulay, certain conditions guarantee that this complex
has zero homology.  Here, we develop {\em Euler--Koszul homology}\/
(Definition~\ref{d:ek}) for the class of {\em toric modules}\/
(Definition~\ref{d:toric}), which are slight generalizations of
$\ZZ^d$-graded modules over the semigroup ring~$\CC[\NN A]$.  For any
toric module~$M$, we show in Theorem~\ref{t:qdeg} that the set of
parameters~$\beta$ for which the Euler--Koszul complex has nonzero
higher homology is precisely the analogue for~$M$ of the exceptional
quasi-degree set~$E_A$
defined above for $M = \CC[\NN A]$.

Having now two cohomology theories, one to recover the local
cohomology quasi-degrees for hypergeometric families, and another to
detect rank-jumping parameters for general holonomic families, we link
them in a central result of this article, Theorem~\ref{t:euler-beta}:
for toric modules, these two theories coincide.  Consequently, we
obtain Theorem~\ref{conj:MM} as the special case $M = \CC[\NN A]$ of
Theorem~\ref{t:MM}, which holds for arbitrary toric modules~$M$.  We
deduce in Corollary~\ref{c:GKZ} the equivalence of the
\mbox{Cohen--Macaulay} condition on~$\CC[\NN A]$ with the
\mbox{absence of rank-jumps in the GKZ hypergeometric
system~$H_A(\beta)$}.

As a final comment, let us note that we avoid the explicit computation
of solutions to hypergeometric systems.  This contrasts with the
previously cited constructions of exceptional parameters, which rely
on combinatorial techniques to produce solutions.  It is for this
reason that these constructions contained the assumption that the
semigroup ring~$\CC[\NN A]$ is graded in the usual $\ZZ$-grading, for
this implies that the corresponding hypergeometric systems are {\em
regular holonomic}\/ and thus have solutions expressible as power
series with logarithms, with all the combinatorial control this
provides.  Our use of homological techniques makes the results in this
article valid in both the regular and non-regular cases.


\subsection*{Acknowledgments}
The equivalence of hypergeometric rank constancy and toric
Cohen--Macaulayness was first conjectured by Bernd Sturmfels, who
divulged it in talks and open-problem sessions.  He was also the first
to mention local cohomology in connection with rank-jumps, although
not in the precise form that appears in this article.  His
encouragement and advice have been invaluable to us throughout this
project.

The occasions when all three authors of this article were in the same
place at the same time have been rare, and all the more appreciated
for that.  We would like to thank the organizers of the Joint
International Meeting of the AMS and RSME in Sevilla, Spain, for
providing us with one such opportunity in the Summer of 2003.  LFM and
UW are especially grateful to Francisco Castro-Jim\'enez, Jose
Mar\'{\i}a Ucha, and Mar\'{\i}a Isabel Hartillo Hermoso, who were our
gracious hosts when we stayed in Sevilla for a week after that
conference.  All three authors also intersected at MSRI.  We thank
this institute for the friendly and stimulating research atmosphere it
provides.

\section{Upper semi-continuity of rank in holonomic families}

\label{sec:upper}

The results in the first part of this article
(Sections~\ref{sec:upper} and~\ref{sec:criteria}) deal with general
modules over the Weyl algebra, without restricting to the
hypergeometric realm.  In this section we define the notion of a {\em
holonomic family}\/ of $D$-modules, and show that the holonomic rank
constitutes an upper semi-continuous function on such a family.

Throughout this article, $\del = \del_1,\dots ,\del_n$ refers to the
partial derivation operators with respect to the variables $x= x_1,
\dots ,x_n$.  Writing $\delta_{ij}$ for the Kronecker delta, so that
$\delta_{ij} = 1$ if $i=j$ and $\delta_{ij} = 0$ otherwise, the Weyl
algebra~$D$ is the quotient of the free associative $\CC$-algebra on
$\{x_i,\del_i\}_{i=1}^n$ by the two-sided ideal
$\ideal{x_ix_j-x_jx_i,\, \del_i\del_j-\del_j\del_i,\, \del_i x_j- x_j
\del_i - \delta_{ij}}_{i,j=1}^n$.

Every left $D$-module is also a module over the commutative
subalgebra~$\CC[x]$ of~$D$.  If~$N$ is any $\CC[x]$-module or sheaf of
$\CC[x]$-modules on some space, let \mbox{$N(x) = \CC(x)
\otimes_{\CC[x]} N$} denote the localization by inverting all
polynomials in the $x$-variables.  We shall also have occasion to
consider modules, sheaves of modules, and schemes defined over~$\CC$,
and if~$N$ is such an object then $N(x)$ denotes its base extension to
the field~$\CC(x)$ of rational functions.

Our focus is on {\em holonomic}\/ $D$-modules, which are by definition
finitely generated left $D$-modules~$N$ such that $\ext^j_D(N,D) = 0$
whenever $j \neq n$.  The holonomic modules form a full subcategory of
the category of $D$-modules that is closed under taking extensions,
submodules, and quotient modules.  When $N$ is a holonomic module, the
$\CC(x)$-vector space~$N(x)$ has finite dimension, and this dimension
equals the {\em holonomic rank}\/~$\rk(N)$ by a celebrated theorem of
Kashiwara, see \cite[Thm.~1.4.19, Cor.~1.4.14]{SST}.  We note that
rank is hence additive in short exact sequences of holonomic modules.

We are interested in families of \mbox{$D$-modules} parameterized by a
Noetherian complex algebraic variety~$B$ with structure sheaf $\calO_B$.
If $\beta \in B$ we denote by $\frakp_\beta$
the prime ideal (sheaf) of~$\beta$, and by~$\kappa_\beta$ the residue
field of the stalk $\calO_{B,\beta}$, so $\kappa_\beta =
\calO_{B,\beta}/\frakp_\beta\calO_{B,\beta}$. 

Consider the sheaf $D \otimes_\CC \calO_B$ of noncommutative
$\calO_B$-algebras on~$B$.  By a {\em coherent}\/ sheaf of left $(D
\otimes_\CC \calO_B)$-modules we mean a quasi-coherent sheaf of
$\calO_B$-modules on~$B$ whose sections over each open affine subset
$U \subseteq B$ are finitely generated over the ring of global
sections of $D \otimes_\CC \calO_U$.  The sheaf $D \otimes_\CC
\calO_B$ contains the subsheaf $\calO_B[x] = \CC[x] \otimes_\CC
\calO_B$ of commutative polynomials over~$\calO_B$, whose localization
at $\ideal{0}\in\spec(\CC[x])$ is by our conventions $\calO_B(x)$.
The sheaf-spectrum of~$\calO_B(x)$ is the base-extended scheme $B(x) =
\spec\CC(x) \times_{\spec\CC} B$, which can be considered as a scheme
over~$\CC(x)$ or as fibered over~$B$.

\begin{defn} \label{d:hol}
A {\em holonomic family}\/ over a complex scheme~$B$ is a coherent
sheaf~$\tilde\calM$ of left $(D \otimes_\CC \calO_B)$-modules whose
fibers $\calM_\beta = \tilde\calM \otimes_{\calO_B} \kappa_\beta$ are
holonomic $D$-modules for all $\beta \in B$, and whose {\em rank
sheaf}\/ $\tilde\calM(x)$ is coherent on the scheme~$B(x)$.
\end{defn}

Note that every holonomic family over~$B$ is generated by its global
sections if $B$ is affine.

\begin{example} \label{ex:CCd}
When $B$ equals the complex vector space~$\CC^d$ of dimension $d$, the
structure sheaf~$\calO_B$ has global sections~$\CC[b]$, the
commutative polynomial ring in variables $b = b_1,\ldots,b_d$.  A
holonomic family $\tilde\calM$ over~$\CC^d$ can be represented by its
global sections, a left $D[b]$-module~$\calM$ with finitely many
generators and relations.  In order to be a holonomic family in our
sense, $\calM_\beta=\calM\otimes_{\calO_B}\kappa_\beta$ has to be a
holonomic $D$-module for all $\beta \in \CC^d$ while the global
section module $\calM(x) = \CC(x) \otimes_{\CC[x]} \calM$ of the rank
sheaf has to be a finitely generated $\CC[b](x)$-module.
\end{example}

The rank sheaf $\tilde\calM(x)$ is a sheaf of $\calO_B(x)$-modules
on~$B(x)$, but ignoring the process of pushing $\tilde\calM(x)$ down
to $B$ we abuse notation and speak of the fiber of~$\tilde\calM(x)$
over $\beta \in B$.  Let $B_\CC\subseteq B$ be the $\CC$-valued closed
points of~$B$ (that is, points with residue field~$\CC$).

\begin{prop} \label{p:rank}
For a holonomic family $\tilde\calM$ over\/ a scheme~$B$, the fiber
of~$\tilde\calM(x)$ over each point $\beta \in B$ is~$\calM_\beta(x)$,
which is a $\CC(x)$-vector space of dimension\/~$\rk(\calM_\beta)$ if
$\beta \in B_\CC$.
\end{prop}
\begin{proof}
Tensoring with the rational functions~$\CC(x)$ over~$\CC[x]$ commutes
with the passage from~$\tilde\calM$ to~$\calM_\beta = \tilde\calM
\otimes_{\calO_B} \kappa_\beta$.
\end{proof}

The next lemma says that Zariski closed subsets of the base-extended
scheme~$B(x)$ descend to Zariski closed subsets of the original
scheme~$B$.  We use the term {\em prime ideal}\/ to refer to the
kernel of the morphism of structure sheaves associated to any map from
the spectrum of a field to~$B$.

\begin{lemma} \label{l:closed}
For a scheme~$B$ defined over\/~$\CC$, the map taking each prime
ideal\/~\mbox{$\pp \subset \calO_B$} to its extension $\pp \calO_B(x)$
constitutes a continuous injection \mbox{$B \to B(x)$} of
topological~spaces.
\end{lemma}
\begin{proof}
Note first that if $\pp \subset \calO_B$ is prime, then the extension
$\pp\calO_B(x)$ is prime.  The lemma is equivalent to saying that, for
every subset \mbox{$Y \subseteq B(x)$} that is Zariski closed, the set
of points $\beta \in B$ whose prime ideal sheaves~$\pp_\beta$
in~$\calO_B$ extend to prime ideals $\pp_\beta \calO_B(x) \in Y$ is
Zariski closed in~$B$.  This statement is local on~$B$, so we may
assume that~$B$ is affine.

Suppose that $Y$ is the variety of a set~$\calF$ of global sections
of~$\calO_B(x)$.  Any free $\CC$-basis for~$\CC(x)$ is also a free
$\calO_B$-basis for $\calO_B(x)$; choose such a basis.  For each
global section $f \in \calF$, let $\calG_f \subseteq \calO_B$ be its set
of nonzero coefficients in the chosen basis.  Then a prime $\pp
\subset \calO_B$ satisfies $\pp\calO_B(x) \in Y$ if and only if $\pp$
contains the sets $\calG_f$ for all $f \in \calF$.%
\end{proof}


\begin{example}
Lemma~\ref{l:closed} does not say that $B$ is closed in~$B(x)$;
indeed, the closure of the image is all of~$B(x)$, whereas there are
always points in \mbox{$B(x) \minus B$}.  Neither does
Lemma~\ref{l:closed} say that the morphism $B(x) \to B$ takes closed
sets to closed sets.  For example, if $n = d = 1$ in
Example~\ref{ex:CCd}, then the variety of the ideal $\<b-x\> \subset
\CC[b](x)$ is a Zariski closed point in~$B(x)$, whereas the image of
this point in~$B = \CC^1$ is the generic point.  Note that the set of
prime ideals $\pp \subset \CC[b]$ whose extensions $\pp \CC[b](x)$
contain $\<b-x\>$ is empty, and therefore Zariski closed.
\end{example}

We now come to the main result of this section.  To state it, we call
$\beta\in B$ {\em rank-jumping}\/ if (i)~$\beta \in B_\CC$ is a
$\CC$-valued point, and (ii)~for all open sets~$U$ containing~$\beta$,
the rank of~$\calM_\beta$ is strictly greater than the minimal
holonomic rank attained by any fiber~$\calM_{\beta'}$ for $\beta'\in
U_\CC$.

\begin{thm} \label{t:semicont}
If~$\calM$ is a holonomic family over a scheme~$B$, then the function
$\beta \mapsto \rk(\calM_\beta)$ is upper semi-continuous on both~$B$
and on~$B_\CC$ (endowed with the subspace topology)\/.
In particular, the locus of rank-jumping points $\beta \in B_\CC$ is
closed in~$B_\CC$.%
\end{thm}
\begin{proof}
The function sending each point $\beta(x) \in B(x)$ to the
$\CC(x)$-dimension of the fiber
$\tilde\calM(x)\otimes\kappa_{\beta(x)}$ is upper-semicontinuous
on~$B(x)$ because $\tilde\calM(x)$ is coherent on~$B(x)$; see
\cite[III.12.7.2]{Hartshorne}.  Therefore, given an integer~$i$, the
subset of points in~$B(x)$ on which this fiber dimension is at
least~$i$ is Zariski closed.  Lemma~\ref{l:closed} shows that the
corresponding subset of points in~$B$ is Zariski closed in~$B$.
Proposition~\ref{p:rank}, which says that the fiber dimensions over
$\beta \in B_\CC$ are holonomic ranks, completes the proof.%
\end{proof}

\begin{remark} \label{rk:drops}
Without the coherence hypothesis on~$\tilde\calM(x)$ over~$B(x)$, the
conclusion of Theorem~\ref{t:semicont} can be false, even if $B$ is of
finite type over~$\CC$ and all of the fibers~$\calM_\beta$ for $\beta
\in B_\CC$ are holonomic.  For an example, consider the setup in
Example~\ref{ex:CCd} with $n = d = 1$, and take $\calM =
D[b]/\<b\del-1\>$.  When $\beta \neq \ideal{b-0}$, the fiber
over~$\beta$ is the rank~$1$ holonomic module corresponding to the
solution~$x^{1/\beta}$.  But when $\beta = \ideal{b-0}$ the fiber
of~$\tilde\calM$ is zero.  Hence the rank actually drops on the closed
subset $\{0\} \subset \CC^1$.  \mbox{See Example~\ref{ex:drops} for
further~details}.
\end{remark}

\begin{remark} \label{rk:sol}
Fix a holonomic family $\tilde \calM$.  The 
semicontinuity of holonomic rank in Theorem~\ref{t:semicont}
suggests that a ``solution sheaf''~$\sol(\tilde\calM)$, constructed as
below, might be an algebraic coherent sheaf on~$B$; a~priori, it
can only be expected to be an analytic sheaf.  See Remark~\ref{rk:MM} for
further comments on this issue, as it relates to hypergeometric
systems.

For the construction, suppose there is a vector $v \in \CC^n$ such
that the singular locus of~$\calM_\beta$ does not contain~$v$ for any
$\beta\in B_\CC$ (this is certainly possible locally on~$B$).  Let the
point~$v$ have maximal ideal $\<x_1-v_1,\dots, x_n-v_n\>$ in~$\CC[x]$.
The $D$-module restriction of~$\tilde\calM_\beta$ to~$v$ is the
derived tensor product over~$D$ with the right module
$D/\<x_1-v_1,\dots, x_n-v_n \>D$.  This restriction is naturally dual
to the space of formal power series solutions of~$\tilde\calM_\beta$
at~$v$.  But as~$v$ is a regular point, formal and convergent solutions
of $\tilde\calM_\beta$ are identical.  Since tensor products commute
in~$x$ and~$b$, restricting the family~$\tilde\calM$ gives rise to an
$\calO_B$-module whose fiber over each point $\beta\in B_\CC$ is
naturally dual to the solution space of~$\tilde\calM_\beta$.  Taking
duals therefore yields the analytic
$\calO_B$-module~$\sol(\tilde\calM)$.
\end{remark}

\section{Rank-jumps as failures of flatness}

\label{sec:criteria}

For applications to hypergeometric systems, we need some concrete
criteria to help us apply the results of the previous section.  First
we characterize rank-jumps in holonomic families over reduced schemes
using Koszul homology in a commutative setting (Corollary~\ref{c:cm}).
Then we provide a criterion for a family of $D$-modules to satisfy the
coherence property required of a holonomic family
(Proposition~\ref{p:filter}); we will appeal to it in
Section~\ref{sec:global}.

In the following statement, we use the standard notion of {\em
reduced}\/ for a scheme to mean that the coordinate rings of its
affine open subschemes have no nilpotent elements.

\begin{thm} \label{t:cm}
Let~$\tilde\calM$ be a holonomic family over a reduced scheme~$B$ of
finite type over\/~$\CC$, and fix a $\CC$-valued closed point $\beta
\in B_\CC$.  The parameter~$\beta$ is rank-jumping if and only if\/
$\Tor_i^{\calO_B}(\tilde\calM(x),\kappa_\beta)$ is nonzero for some
homological degree $i > 0$.
\end{thm}
\begin{proof}
Since $B$ is of finite type, the set $B_\CC$ is dense in~$B$ via the
inclusion of Lemma~\ref{l:closed}.  As $B$ is always dense in~$B(x)$,
we deduce that $B_\CC$ is also dense in~$B(x)$.  Therefore
\begin{stareqn} \label{r}
  \min_{\beta\in B_\CC}\{\rk(\calM_\beta)\} &=& \min_{\beta\in
  B(x)}\{\dim_{\CC(x)}\tilde\calM_\beta(x)\}
\end{stareqn}%
by Proposition~\ref{p:rank} and the upper-semicontinuity of fiber
dimension \cite[III.12.7.2]{Hartshorne}.  By
\cite[Exercise~II.5.8]{Hartshorne}, the points $\beta(x) \in B(x)$
over which the fiber dimension equals the quantity given in~(\ref{r})
coincide with the points~$\beta(x)$ near which~$\tilde\calM(x)$ is
locally free.  But local freeness and flatness agree for coherent
sheaves by \cite[Proposition~III.9.2]{Hartshorne}.  Flatness in turn
is characterized by the vanishing of
$\Tor_i^{\phantom.}{}^{\hspace{-1ex}\calO_{B(x)}}
(\tilde\calM(x),\kappa_{\beta(x)})$ for all \mbox{$i > 0$}
\cite[Theorem~6.8]{Eisenbud}, and for $\beta \in B_\CC$,
$\Tor_i^{\phantom.}{}^{\hspace{-1ex}\calO_{B(x)}}
(\tilde\calM(x),\kappa_{\beta(x)})\cong
\Tor_i^{\calO_B}(\tilde\calM(x),\kappa_\beta)$.%
\end{proof}

\begin{remark}
The result of Theorem~\ref{t:cm} is false if $B$ is not reduced.
Consider for example $B=\spec(\CC[b]/\ideal{b^2})$, and let
$\tilde\calM$ be the sheaf induced by the module
$D[b]/\ideal{\del,b}$. Since there is only one point in $B$, there is
no possibility of a rank-jump. On the other hand,
$\tor_i^{\calO_B}(\tilde \calM(x),\kappa_\beta)=D(x)/\ideal{\del}$ for
all $i$.
\end{remark}

Theorem~\ref{t:cm} has an interpretation via Koszul homology.  For
notation, suppose $T$ is a commutative ring.  For a sequence $\boldy =
y_1,\ldots,y_d$ of elements in~$T$, we write the Koszul complex
$K_\spot(\boldy)$ \cite[Chapter~1.6]{BH93} with lowered indices
decreasing from~$d$ to~$0$.  For any $T$-module~$N$ set
$K_\spot(\boldy;N) = K_\spot(\boldy) \otimes_T N$, and abbreviate the
$i^\th$ Koszul homology of~$N$ as $H_i(\boldy;N) = H_i
(K_\spot(\boldy;N))$.  Call $\boldy$ a {\em regular sequence}\/ in~$T$
if $T/\ideal{\boldy}$ is a nonzero module and $y_i$ is a
non-zerodivisor on $T/\<y_1,\ldots,y_{i-1}\>$ for $i = 1,\ldots,d$.
Koszul cohomology for regular sequences gives Tor groups
\cite[Exercise~17.10]{Eisenbud}, so we get the following.

\begin{cor} \label{c:cm}
Let~$\calM$ be the global sections of a holonomic family over a
reduced affine variety~$B$ with finitely generated coordinate
ring\/~$\CC[B]$.  If~a regular sequence~$\boldy$~in $\CC[B]$ generates
the ideal of a $\CC$-valued closed point $\beta \in B_\CC$, then
$\beta$ is rank-jumping if and only if the Koszul complex
$K_\spot(\boldy;\calM(x))$ has nonzero homology $H_i(\boldy;\calM(x))$
for some $i > 0$.\qed
\end{cor}

Note that $K_\spot(\boldy;\calM(x))$ can be obtained either by
tensoring $K_\spot(\boldy)$ with $\calM(x)$ over~$\CC[B]$, or by
viewing $\boldy$ as a sequence in~$\CC[B](x)$ and tensoring
over~$\CC[B](x)$.

%
\smallskip

We now turn to the coherence criterion.  We refer to
\cite{Bjork-Diffops,SST} for more information about filtrations on
$D$-modules and their associated graded objects.  On $D$ define the
{\em order filtration}\/ by taking its $k^\th$ level to be the vector
space of all expressions $\sum_{\nu} p_{\nu}(x) \del^\nu$ in which
the monomials~$\del^\nu$ for $\nu \in \NN^n$ have total
degree~\mbox{$|\nu| \leq k$}.  Note that each filtration level is a
finitely generated $\CC[x]$-module.  The associated graded object is
the commutative polynomial ring $\gr (D) = \CC[x,\xi]$ in the
variables $x = x_1,\ldots,x_n$ and $\xi = \xi_1,\ldots,\xi_n$.

The order filtration extends to any free module~$D^r$ with a fixed
basis, and to left submodules $K \subseteq D^r$, by letting the
$k^\th$ level be spanned by elements whose coefficients in the given
basis have total degree at most~$k$ in~$\del$.  Given a presentation
$N = D^r/K$, the left module~$N$ has an induced filtration with
associated graded module $\gr(N) = \gr(D^r)/\gr(K)$.  This naturally
graded $\gr(D)$-module depends on the choice of presentation.

Note that since $\CC(x)$ is flat over~$\CC[x]$, the module~$\gr(N)(x)$
is isomorphic to~$\gr(N(x))$ if we extend the filtration above in the
obvious way to $N(x)=(D^r/K)\otimes_{\CC[x]}\CC(x)$.  It is a
fundamental fact of $D$-module theory, going back to Kashiwara, that
for a holonomic $D$-module~$N$, the number $\dim_{\CC(x)}(\gr(N)(x))$
is independent of the presentation of~$N$ and equals the holonomic
rank of~$N$ \cite[Theorem~1.4.19 and Definition~1.4.8]{SST}.

Our primary use of order filtrations will be on {\em families}\/ of
$D$-modules over an affine scheme~$B$---that is, on left modules over
the ring $D[B] = D \otimes_\CC \CC[B]$.  The order filtration
generalizes naturally to a filtration on $D[B]$ with the property that
each level is a finitely generated $\CC[B][x]$-module.  The associated
graded ring of~$D[B]$ is the polynomial ring $\CC[B][x,\xi]$ in $2n$
variables over the coordinate ring~$\CC[B]$.  A~choice of presentation
for a $D[B]$-module~$\calM$ determines an associated graded
$\CC[B][x,\xi]$-module~$\gr(\calM)$.  The formation of such associated
graded structures commutes with the tensor product with the flat
$\CC[x]$-module $\CC(x)$.

\begin{prop} \label{p:filter}
For an affine scheme~$B$ and a finitely generated\/
\mbox{$D[B]$-module}~$\calM$, the module~$\calM(x)$ is finitely generated
over\/~$\CC[B](x)$ provided that $(\gr(\calM))(x)$~is.
\end{prop}
\begin{proof}
As $B$ is affine, the modules $\calM$, $\calM(x)$ and $\gr(\calM(x))$
are generated by global sections.  Any lift of set of generators for
the $\CC[B](x)$-module $(\gr\calM)(x)$ is a set of generators
for~$\calM(x)$, because the order filtration of~$\calM(x)$ does not
descend infinitely.%
\end{proof}


We will apply the above result when all fibers~$\calM_\beta$ for
\mbox{$\beta \in B_\CC$} are holonomic $D$-modules to conclude that
$\calM$ constitutes (the global sections of) a~holonomic family
over~$B$.  It~should be pointed out that the finiteness condition
on~$\calM$ is necessary: even if the fibers~$\calM_\beta$ are
holonomic, the fiber $(\gr\calM)(x)_\beta$ of the associated graded
module~$(\gr\calM)(x)$ over \mbox{$\beta \in B_\CC$} is a
$\CC(x)$-vector space whose dimension need not be equal
to~$\rk(\calM_\beta)$.

\begin{example} \label{ex:drops}
Continue with~$\calM = D[b]/\< b\del -1\>$ as in
Remark~\ref{rk:drops}.  Let $\calM(x)_0 \subseteq \calM(x)_1 \subseteq
\cdots$ denote the order filtration of~$\calM(x)$.  Then $\calM(x)_1 =
\calM(x)$ locally near every parameter \mbox{$\beta \in \CC$} except
for $\beta = 0$.  At~\mbox{$\beta = 0$}, in contrast,
$\calM(x)_k/\calM(x)_{k-1}$ is minimally generated by $\xi^k$, even
though $\rk(\calM_0) = 0$.  In~fact, $(\gr\calM)(x)$ is
the direct sum of the rank~$1$ free module~$\CC[b](x)$ with
$\bigoplus_{k \geq 1} \CC[b](x)/\<b\>$.
\end{example}

\section{Euler--Koszul homology of toric modules}

\label{sec:ek}

In this section we introduce generalized $A$-hypergeometric systems in
the sense indicated in the introduction.  After reviewing some basic
facts of GKZ hypergeometric systems, we provide foundations for
Euler--Koszul homology of what we call {\em toric $\ZZ^d$-graded
modules}.

For the rest of this paper, fix a $d\times n$ integer
matrix~$A=(a_{ij})$ of rank~$d$.  We emphasize that we do not assume
that the columns $a_1,\ldots,a_n$ of~$A$ lie in an affine hyperplane.
However, we do assume that $A$ is {\em pointed}, meaning that
$a_1,\ldots,a_n$ lie in a single open half-space of~$\RR^d$. This
guarantees that the semigroup
\begin{eqnarray*}
  \NN A &=& \Big\{\sum_{i=1}^n \gamma_i a_i \mid
  \gamma_1,\ldots,\gamma_n\in\NN\Big\}
\end{eqnarray*}
has no units.  (Pointedness will come into play in the proof of
Theorem~\ref{t:qdeg}, where it is used to ensure that local duality
holds.)  The {\em semigroup ring}\/ associated to the $d \times n$
matrix~$A$ is $\SA = \CC[\NN A] \cong R/I_A$, where $R =
\CC[\del_1,\ldots,\del_n]$~and
\begin{eqnarray*}
  I_A &=& \<\del^\mu-\del^\nu \mid \mu,\nu\in\ZZ^n, A\cdot \mu =
  A\cdot \nu\>
\end{eqnarray*}
is the {\em toric ideal}\/ of~$A$.  Notice that $\SA$ and $R$ are
naturally graded by~$\ZZ^d$ if we define $\deg(\del_j) = -a_j$, the
negative of the $j^\th$ column of~$A$.

Our choice of signs in the $\ZZ^d$-grading of~$R$ is compatible with a
$\ZZ^d$-grading on the Weyl algebra~$D$ in which $\deg(x_j) = a_j$ and
\mbox{$\deg(\del_j) = -a_j$}.  Under this $\ZZ^d$-grading, the $i^\th$
{\em Euler operator}\/ $E_i = \sum_{j=1}^n a_{ij} x_j\del_j\in D$ is
homogeneous of degree~$0$ for $i = 1,\ldots,d$.  Given a vector
\mbox{$\beta \in \CC^d$}, we write $E-\beta$ for the sequence
$E_1-\beta_1,\ldots,E_d-\beta_d$.

\begin{defn}\label{def:GGKZ}
The {\em $A$-hypergeometric GKZ system}\/ with {\em
parameter~$\beta$}\/ is the left ideal
\begin{eqnarray*}
  H_A(\beta) &=&  D\cdot\ideal{I_A,E-\beta}
\end{eqnarray*}
in the Weyl algebra~$D$.  The {\em $A$-hypergeometric $D$-module}\/
with {\em parameter~$\beta$}\/ is
\begin{eqnarray*}
  \calM^A_\beta &=& D/\!H_A(\beta).
\end{eqnarray*}  
\end{defn}

Results of \cite{GKZ89,Adolphson-Duke94,Hotta,SST} imply that the
module~$\calM^A_\beta$ is holonomic of nonzero rank.
$A$-hypergeometric systems constitute an important class of
$D$-modules, playing a role similar to that of toric varieties in
algebraic geometry, since they possess enough combinatorial
underpinning to make calculations feasible, and enough diversity of
behavior to make them interesting as a test class for conjectures and
computer \mbox{experimentation}.
%

If~\mbox{$y \in N_\alpha$} is homogeneous of degree~$\alpha$ in a
$\ZZ^d$-graded $D$-module~$N$, write $\deg_i(y) = \alpha_i$.  In
particular, for any homogeneous element $P \in D$,
\begin{eqnarray*}
  E_iP-PE_i &=& \deg_i(P)P.
\end{eqnarray*}
The displayed equation shows that if $N$ is any $\ZZ^d$-graded left
$D$-module, then the map $E_i - \beta_i: N \to N$ sending the
homogeneous element $y \in N$ to
\begin{eqnarray*}
  (E_i - \beta_i) \circ y &=& (E_i-\beta_i-\deg_i(y)) y
\end{eqnarray*}
determines a $D$-linear endomorphism of~$N$ when extended
$\CC$-linearly to inhomogeneous elements of~$N$.  This endomorphism is
functorial, in the sense that it commutes with degree~$0$
homomorphisms $N \to N'$ of $\ZZ^d$-graded left $D$-modules.
Moreover, the endomorphisms for the various Euler operators commute:
$[E_i-\beta_i,E_j-\beta_j]$ is the zero endomorphism on any
$\ZZ^d$-graded left $D$-module for all $i,j = 1,\ldots,d$ and all
complex numbers $\beta_i, \beta_j$.

\begin{defn} \label{d:ek}
Fix $\beta\in\CC^d$ and a $\ZZ^d$-graded $R$-module~$N$.  The {\em
Euler--Koszul complex}\/ $\calK_\spot(E-\beta;N)$ is the Koszul
complex of left $D$-modules defined by the sequence $E-\beta$ of
commuting endomorphisms on the left $D$-module $D\otimes_R N$
concentrated in homological degrees $d$ to $0$.  The $i^\th$ {\em
Euler--Koszul homology}\/ of~$N$ is
$\calH_i(E-\beta;N)=H_i(\calK_\spot(E-\beta;N))$.
\end{defn}

This complex was used before in the special case $N=\SA$ in
\cite{GKZ89, Adolphson-Duke94,Adolphson-Rend} in order to
study the solutions to~$H_A(\beta)$ for special parameters~$\beta$.
We use the script~`$\calK$' for Euler--Koszul complexes instead of the
usual~`$K$' for ordinary Koszul complexes to emphasize that the maps
in $\calK_\spot(E-\beta;N)$ are homomorphisms of $D$-modules. 
The action of an endomorphism in the sequence $E-\beta$ on a
  homogeneous element depends on the degree of that element, in
  contrast to maps in Koszul complexes over commutative graded
  rings.
Nonetheless, Euler--Koszul complexes behave much like ordinary Koszul
complexes.  To see why, let
$\theta_i=x_i\del_i\in D$.  In the $\ZZ^d$-graded commutative subalgebra
$\Theta=\CC[\theta_1,\ldots,\theta_n]$ of~$D$
every element has degree \mbox{$0 \in \ZZ^d$}.
Consequently, $N_\alpha$ is a left $\Theta$-module whenever $N$ is a
$\ZZ^d$-graded $D$-module. We hence obtain the following.

\begin{lemma} \label{l:honest}
Let $N$ be a $\ZZ^d$-graded $R$-module and $\alpha \in \ZZ^d$.  The
$\ZZ^d$-graded degree~$\alpha$ part $\calK_\spot(E-\beta;N)_\alpha$ of
the Euler--Koszul complex agrees with the ordinary Koszul complex
$K_\spot(E-\beta-\alpha; (D \otimes_R N)_\alpha)$ constructed over the
polynomial subalgebra $\Theta \subseteq D$.  That is,
\begin{eqnarray*}
  \calK_\spot(E-\beta;N) &=&
  \bigoplus_{\alpha\in\ZZ^d}K_\spot(E-\beta-\alpha;(D\otimes_R
  N)_\alpha),
\end{eqnarray*}
the right hand side being a direct sum of ordinary Koszul complexes of
$\Theta$-modules.\qed
\end{lemma}

Since $D$ is a free $R$-module, the Euler--Koszul complex constitutes
an exact functor from the category of $\ZZ^d$-graded $R$-modules with
$\ZZ^d$-graded morphisms of degree $0$ to the category of bounded
complexes of $\ZZ^n$-graded left $D$-modules with $\ZZ^d$-graded
morphisms of degree $0$.  Consequently, short exact sequences of
$\ZZ^d$-graded modules over~$R$ give rise to long exact sequences of
Euler--Koszul homology.

We now specify a subcategory of the category of $R$-modules that is
central to what follows.  This subcategory will contain the semigroup
ring~$\SA$ as well as its monomial (that is, $\ZZ^d$-graded) ideals,
and it will be closed under taking $\ZZ^d$-graded submodules, quotient
modules, and extensions.  In the following section we then study the
Euler--Koszul homology of objects in this category and establish
strong (non)vanishing conditions.

\begin{notn} \label{n:toric}
By~a {\em face of~$A$} we mean a set~$F$ of columns of $A$ minimizing
some linear functional on~$\NN A$. We interpret $F$ as a submatrix of
$A$.  The {\em dimension}\/ $\dim(F)$ of a face~$F$ equals the rank of
the subgroup $\ZZ F \subseteq \ZZ^d$ it generates.  Let~$\SF$ be the
semigroup ring generated by $F$.  Thus $\SF \cong R/I_F$, where
\mbox{$I_F = \IAF + \<\del_j \mid a_j \not \in F\>$} is the prime
ideal obtained by starting with the toric ideal~$\IAF$ (defined as the
kernel of the surjection \mbox{$\CC[\del_j \mid a_j \in F] \to \SF$})
and adding the variables corresponding to the columns of~$A$
outside~of~$F$.
\end{notn}

\begin{defn} \label{d:toric}
A $\ZZ^d$-graded $R$-module~$M$ is {\em toric}\/ if it has a
{\em toric filtration}
$$
  0 = M_0 \subset M_1 \subset \cdots \subset M_{\ell-1} \subset M_\ell
  = M,
$$
meaning that $M_k/M_{k-1}$ is, for each~$k$, a $\ZZ^d$-graded
translate of $\SFk$ for some face~$F_k$ of~$\NN A$.  We say that $M$
has {\em toric length}\/~$\ell$ if the minimal length of a toric
filtration for $M$ is~$\ell$.
\end{defn}

\begin{remark} \label{rem:toric}
Most toric modules have many different toric filtrations, and usually
more than one of these has minimal length.  A toric module of toric
length~$1$ is simply a $\ZZ^d$-graded translate $\SF(-\alpha)$
of~$\SF$ generated in degree~\mbox{$\alpha \in \ZZ^d$}, for some
face~$F$ of~$A$.
\end{remark}

\begin{example} \label{ex:toric}
If $N$ is a finitely generated $\ZZ^d$-graded $\SF$-module for some
face~$F$ of~$A$, then $N$ is toric.  To see this, argue by induction
on~$\dim(N)$, the Artinian case $\dim(N) = 0$ being trivial.  If
$\dim(N)>0$, then $N$ has a submodule of the form $\SFp(-\alpha)$ with
$F'\subseteq F$ being some face of $A$ of dimension $\dim(N)$. If we
replace $N$ by $N/\SFp(\alpha)$ then a finite number of such steps
will reduce the dimension of $N$ by at least one.
\end{example}

Example~\ref{ex:toric} implies that $\ZZ^d$-graded submodules,
quotients, and extensions of toric modules are toric, since they have
a composition chain by shifts of modules of the form~$\SF$.


If $F$ is a face of~$A$, then we write $E^F_i$ for the Euler operator
obtained by erasing all terms $a_{ij}x_j\del_j$ from~$E_i$ such that
$a_j \not\in F$.
We will now consider the actions of Euler endomorphisms $E_i-\beta_i$
on $D \otimes_{R} M$ for an $\SF$-module~$M$.  To describe such
actions, let~$D_F$ and~$D_\ol F$ be the Weyl algebras on the variables
$\{x_F,\del_F\}$ and $\{x_\ol F,\del_\ol F\}$ corresponding to the
columns~$a_j \in F$ and $a_j \not\in F$, respectively.  Then for each
$\SF$-module $M$ we have
\begin{stareqn} \label{iso}
  D \otimes_{R} M &\cong& \CC[x_\ol F] \otimes_\CC (D_F
  \otimes_{\CC[\del_F]} M)
\end{stareqn}%
as left modules over $D = D_\ol F \otimes_\CC D_F$, identifying
$x_F^\mu x_{\ol F}^{\ol{\mu}} \del_F^\nu\del_{\ol F}^{\ol{\nu}}\otimes
m$ with $x_{\ol F}^{\ol{\mu}}\otimes (x_F^\mu\del_F^\nu\otimes
\del_{\ol F}^{\ol{\nu}}m)$, where of course both sides are zero if
$|\ol\nu|>0$.

Suppose that $M$ is a $\ZZ^d$-graded $\SF$-module and interpret it as
a $\CC[\del_F]$-module (rather than as a module over a quotient of
$R$). Note that $E^F_i$
induces an Euler endomorphism on $D_F\otimes_{\CC[\del_F]} M$ by
setting $E^F_i\circ (P\otimes m)=\sum_{j\in
F}(a_{ij}x_j\del_j-\deg_i(P\otimes m))P\otimes m$ for all homogeneous
elements $P\in D_F$ and $m\in M$.

\begin{lemma} \label{l:EF}
If~$M$ is a finitely generated $\ZZ^d$-graded $\SF$-module, then under
the isomorphism~(\ref{iso}) above, the action of the endomorphism
$E_i-\beta_i$ on $D \otimes_{R} M$ coincides with the action of the
endomorphism $E^F_i-\beta_i$ on the right hand factor of\/ $\CC[x_\ol
F] \otimes_\CC (D_F \otimes_{\CC[\del_F]} M)$:
\begin{eqnarray*}
(E_i-\beta_i)\circ (x_F^\mu x_{\ol F}^{\ol{\mu}} \del_F^\nu\del_{\ol
  F}^{\ol{\nu}}\otimes m) &=& x_{\ol F}^{\ol{\mu}}\otimes
  (E^F_i-\beta_i)\circ (x_F^\mu\del_F^\nu\otimes \del_{\ol
  F}^{\ol{\nu}}m).
\end{eqnarray*}
\end{lemma}
\begin{proof}
Let $P=x_F^\mu x_{\ol F}^{\ol{\mu}} \del_F^\nu\del_{\ol F}^{\ol{\nu}}$
and $P'=x_F^\mu \del_F^\nu\del_{\ol F}^{\ol{\nu}}$.  Note that
$[E_i^F,x_{\ol F}^{\ol{\mu}}] =0$ and therefore $[E_i-E_i^F,x_{\ol
F}^{\ol{\mu}}]=\deg_i(x_{\ol F}^{\ol{\mu}}) x_{\ol F}^{\ol{\mu}}$.  We
then compute
\begin{eqnarray*}
(E_i-\beta_i)\circ (P\otimes m)
&=&(E_i-\beta_i-\deg_i(P\otimes m))\cdot P\otimes m
\\
&=&\bigl((E_i-E_i^F)-\deg_i(x_{\ol F}^{\ol{\mu}})\bigr)P\otimes m
\\
& &\phantom{x}+\bigl(E_i^F-\beta_i-\deg_i(P'\otimes m)\bigr)P\otimes m
\\
&=&0+x_{\ol F}^{\ol\mu}\bigl(E_i^F-\beta_i-\deg_i(P'\otimes m)\bigr)P'\otimes m
\end{eqnarray*}
which is identified with $x_{\ol F}^{\ol{\mu}}\otimes
(E_i^F-\beta_i-\deg_i(P'\otimes m))P'\otimes m$.%
\end{proof}

In analogy with the construction of the $A$-hypergeometric
system~$\calM^A_\beta$ from~$A$, the faces of~$A$ give rise to
holonomic $D$-modules, but sometimes these modules can be zero.

\begin{lemma} \label{l:AF}
Fix a vector $\beta \in \CC^d$ and a face~$F$ of\/~$A$.  The following
are equivalent for the quotient $\calM^F_\beta =
D/\ideal{I_F,E-\beta}$.
\begin{numbered}
\item
The module $\calM^F_\beta$ is nonzero.

\item
$\calM^F_\beta$ is a holonomic $D$-module of nonzero rank.

\item
The parameter $\beta$ (or, equivalently, $-\beta$) lies in the span
(over\/~$\CC$) of the columns of~$F$.
\end{numbered}
\end{lemma}
\begin{proof}
If $\beta$ lies outside the column-span of~$F$, then the $\CC$-linear
span of $E-\beta$ and the products $\{x_j\del_j \mid a_j \not\in F\}
\subset D\cdot I_F$ contains a nonzero scalar, so $\calM^F_\beta$ is zero,
as is its rank.

Suppose now that $\beta$ lies in the column-span of~$F$.  Let $A'_F$
be a submatrix of $A$ composed of $\dim(F)$ many linearly independent
rows of $A_F$. Then $\MAFB = D_F/\ideal{\IAF, E^F-\beta}$ equals
$\MAFBp$.  Hence $\MAFB$ is a holonomic $D_F$-module of nonzero rank,
by \cite{GKZ89,SST}.  Since $\CC[x_\ol F] \otimes_\CC \MAFB \cong
\calM^F_\beta$ by Lemma~\ref{l:EF}, $\calM^F_\beta$ is holonomic over
$D = D_\ol F \otimes_\CC D_F$ and has nonzero rank.%
\end{proof}

We close this section by defining the main object of study of the
remainder of the paper.

\begin{defn}
Let $M$ be a toric $R$-module and~$\beta \in \CC^d$.  The {\em
generalized hypergeometric system}\/ associated to~$M$ and~$\beta$ is
the zeroth Euler--Koszul homology $\calH_0(E-\beta;M)$.
\end{defn}

If $M=\SA$ then $\calH_0(E-\beta;M)$ is the holonomic GKZ-module
$\calM^A_{\beta}$ from Definition~\ref{def:GGKZ}.

\section{Rigidity and holonomicity of Euler--Koszul homology}

We begin our treatment of generalized hypergeometric systems by
showing that they are holonomic, as are all of the higher
Euler--Koszul homology modules of toric modules.

\begin{prop} \label{p:holonomic}
The generalized hypergeometric system $\calH_0(E-\beta;M)$ is
holonomic for every toric module~$M$ and every parameter $\beta \in
\CC^d$.  Consequently, $\calH_i(E-\beta;M)$ is holonomic for all $i>0$
as well.
\end{prop}
\begin{proof}
Consider first the case of $\calH_0(E-\beta;M)$.  The proof is by
induction on the toric length~$\ell$ of~$M$.  When $M = \SF$ for some
face~$F$, use Lemma~\ref{l:AF}.  The general $\ell=1$ case follows
because $\calH_0(E-\beta;\SF(\alpha)) \cong
\calH_0(E-\beta+\alpha;\SF)(\alpha)$.  For $\ell > 1$, the
Euler--Koszul functor applied to a toric short exact sequence $M_{1}
\into M \onto M/M_{1}$ induces a sequence
$$
  \calH_0(E-\beta;M_{1})\to \calH_0(E-\beta;M)\to
  \calH_0(E-\beta;M/M_{1})
$$
that is exact and has holonomic modules at both ends by induction.

Now consider $i>0$.  Let $m\in\ker(\calK_i(E-\beta;M)\to
\calK_{i-1}(E-\beta;M))$, so its coset~$\bar m$ is an element of
$\calH_i(E-\beta;M)$.  Since the Euler--Koszul complex is
$\ZZ^d$-graded, we may assume that $m$ is homogeneous of degree
$\alpha\in\ZZ^d$.  Note that $\calH_i(E-\beta;M)$ is generated by
finitely many such~$\bar m$, since $\calK_i(E-\beta;M)$ is a direct sum of
${\binom{d}{i}}$ copies of $\left(D\otimes_R M\right)$, which is a
Noetherian $D$-module.  It is hence sufficient to prove that $D\cdot
\bar m$ is holonomic.

Consider the Koszul complex $K_\spot(E-\beta-\alpha,(D\otimes_R
M)_\alpha)$ on the (left) $\Theta$-module \mbox{$(D\otimes_R
M)_\alpha$}.  Within this complex, $m$ descends to an element~$\bar m$
of the $i^\th$ homology module $H_i(E-\beta-\alpha,(D\otimes_R
M)_\alpha)$.  We remark that the two usages of $\bar m$ agree in the
sense of the isomorphism in Lemma \ref{l:honest}.  Since $\Theta$ is
commutative, left multiplication by $E_j-\beta_j-\alpha_j$ annihilates
$\bar m$ by \cite[Proposition~1.6.5]{BH93}.  On the other hand, $m$ lies in
the direct sum of ${\binom{d}{i}}$ copies of $(D\otimes_R M)$. Since
$I_A\cdot M=0$, there is some integer $k\in \NN$ such that $(I_A)^k m
= 0$.  It follows that the $D$-module $D\cdot \bar m$ is a
$\ZZ^d$-graded quotient of $D/\ideal{E-\beta-\alpha,(I_A)^k}$. In
particular, $D\cdot \bar m$ has a finite composition series such that
each composition factor is a quotient of
$D/\ideal{E-\beta-\alpha,I_A}\cong \calH_0(E-\beta-\alpha;S_A)$ and
hence is holonomic by the first sentence of the proposition.%
\end{proof}

\begin{defn}
Let~$N$ be any $\ZZ^d$-graded $R$-module.  A~vector $\alpha \in \ZZ^d$
is a {\em true degree}\/ of~$N$, written $\alpha\in\tdeg(N)$, if the
graded piece~$N_\alpha$ is nonzero.  A~vector $\alpha \in \CC^d$ is a
{\em quasi-degree}\/ of~$N$, written $\alpha\in\qdeg(N)$, if $\alpha$
lies in the complex Zariski closure $\qdeg(N)$ of the true degrees
of~$N$ via the natural embedding $\ZZ^d\into\CC^d$.
\end{defn}

We now prove a rigidity property of the Euler--Koszul complex.

\begin{prop} \label{p:HH}
For a toric $R$-module~$M$ and $\beta \in \CC^d$ the following
are equivalent.
\begin{numbered}
\item
$\calH_0(E-\beta;M)$ has holonomic rank~$0$.

\item
$\calH_0(E-\beta;M) = 0$.

\item
$\calH_i(E-\beta;M) = 0$ for all $i \geq 0$.

\item
$-\beta \not\in \qdeg(M)$.
\end{numbered}
\end{prop}
\begin{proof}
2 $\iff$ 3: The failure of the last condition in Lemma~\ref{l:AF} is
equivalent to the $\CC$-span of $E-\beta$ containing an endomorphism
whose action on $D/D\cdot I_F$ is multiplication by a nonzero scalar.  This
observation, together with the isomorphism
$\calH_0(E-\beta;\SF(\alpha)) \cong
\calH_0(E-\beta+\alpha;\SF)(\alpha)$, shows that 2~$\iff$~3 when~$M$
has toric length~$1$.  Now suppose that $M$ has toric
length~\mbox{$\ell > 1$}, and assume the result for modules of smaller
toric length.  For any toric filtration of~$M$, the long exact
sequence of Euler--Koszul homology for $M_1 \into M \onto M/M_1$ shows
that $\calH_0(E-\beta;M)$ surjects onto $\calH_0(E-\beta;M/M_1)$.
Assuming that the former is zero, so is the latter.  By induction, we
find that \mbox{$\calH_i(E-\beta;M/M_1)$} vanishes for all $i \geq 0$,
and hence that $\calH_i(E-\beta;M) \cong \calH_i(E-\beta;M_1)$ for all
$i \geq 0$.  The result now holds for~$M$ because it holds by
induction for~$M_1$.

1 $\iff$ 2: We need that $\calH_0(E-\beta;M) \neq 0$ implies that its
rank is nonzero.  Again use induction on the toric length~$\ell$.
For~\mbox{$\ell=1$} this is Lemma~\ref{l:AF} plus the isomorphism
\mbox{$\calH_0(E-\beta;\SF(\alpha)) \cong
\calH_0(E-\beta+\alpha;\SF)(\alpha)$}, so assume~\mbox{$\ell > 1$}.
Make again the observation that \mbox{$\calH_0(E-\beta;M)$} surjects
onto $\calH_0(E-\beta;M/M_{1})$.  If this latter module is nonzero,
then it has nonzero rank by the \mbox{$\ell = 1$} case, so
$\calH_0(E-\beta;M)$ has nonzero rank, too.  Hence we assume that
$\calH_0(E-\beta;M/M_{1}) = 0$.  Now, using the
equivalence~2~$\iff$~3, we find that $\calH_0(E-\beta;M) \cong
\calH_0(E-\beta;M_{1})$, so again we are done by induction.

2 $\iff$ 4: For the semigroup ring~$\SF$ of a face,
$\calH_0(E-\beta;\SF) = \calM^F_\beta$ is the $D$-module from
Lemma~\ref{l:AF}.  Therefore 2~$\iff$~4 for $M = \SF$.  For
$\ZZ^d$-graded translates of~$\SF$,
$$
\begin{array}{rcccl}
  \left[\calH_0(E-\beta;\SF) \neq 0 \right]&
  \iff&\left[ -\beta \in \qdeg(\SF)\right] &
  \iff&\left[ -\beta-\alpha \in \qdeg(\SF(\alpha))\right]
\end{array}
$$
for $\alpha \in \ZZ^d$, but also $\left[\calH_0(E-\beta;\SF) \neq
0\right] \iff \left[\calH_0(E-\beta-\alpha;\SF(\alpha)) \neq 0\right]$
by definition of the endomorphisms $E_i-\beta_i$.  This proves
2~$\iff$~4 when $M$ has toric length~$1$.  To treat modules of toric
length $\ell > 1$, consider a toric sequence $M_1 \into M \onto
M/M_1$.  Since $\qdeg(M) = \qdeg(M_1) \cup \qdeg(M/M_1)$, we find that
$-\beta \in \qdeg(M)$ if and only if
$$
  \left[-\beta \in \qdeg(M_1)\right]
  \quad \hbox{ or } \quad 
  \left[-\beta \in \qdeg(M/M_1)\right].
$$
This condition is equivalent by induction on~$\ell$ to
$$
  \left[\calH_0(E-\beta;M_1) \neq 0 \right]
  \quad\hbox{ or } \quad
  \left[\calH_0(E-\beta;M/M_1) \neq 0\right]. 
$$
Now $\calH_0(E-\beta;M)$ always surjects onto $\calH_0(E-\beta;M/M_1)$,
and if $\calH_0(E-\beta;M/M_1)$ vanishes then by the equivalence
2~$\iff$~3, $\calH_0(E-\beta;M_1) \cong \calH_0(E-\beta;M)$.
Therefore, this last displayed condition is equivalent to
$\calH_0(E-\beta;M) \neq 0$.%
\end{proof}

\section{Euler--Koszul homology detects local cohomology}

\label{sec:localcoh}

In this section we describe the set of parameters $\beta\in \CC^d$ for
which the Euler--Koszul complex \mbox{$\calK_\spot(E-\beta;N)$} has
nonzero higher homology.  Namely, we identify this set with the
quasi-degrees of the local cohomology of the toric module~$N$.  Our
proof of this result in Theorem~\ref{t:qdeg} uses a spectral sequence,
constructed in Theorem \ref{thm:tot}, that arises from the holonomic
duality functor.  For background on holonomic \mbox{duality we refer
to \cite{Bjork-Diffops}}.

Duality $\DD$ for (complexes of) $D$-modules is the combination of the
derived functor $\RR\!\hom_D(\blank,D)$ of homomorphisms to~$D$
followed by the involution $\tau$ taking~$x^\mu\del^\nu$~to
\begin{eqnarray*}
  \tau(x^\mu\del^\nu) &=& (-\del)^\nu x^\mu.
\end{eqnarray*}
Hence, in order to compute the dual $\DD(N)$ of the module $N$ placed
in homological degree zero, let $F_\spot$ be a $D$-free resolution of
$N$ and apply~$\tau$ to $\hom_D(F_\spot,D)$.  Holonomicity of $N$ is
equivalent to $\DD(N)$ being exact in all cohomological degrees but
$n$.  If $N$ is a module in {\em homological}\/ degree~$k$, then the
{\em co}\/homology $H^i\DD(N)$ is nonzero only for $i = k+n$.
Therefore $\DD$ constitutes an exact functor from the category of
complexes of $D$-modules with holonomic homology to itself, and
$\DD(\DD(N))=N$ for all~$N$.

The rank of a holonomic module~$N$ equals that of its dual: $\rk(N) =
\rk(\DD(N))$.  To be more precise, locally near a nonsingular point of
the complex analytic manifold $\CC^n_{an}$, $N$ is a connection:
$N_{an} = {D_{an}}^r/\ideal{\del_j - C_j \mid j = 1,\ldots,n}$ where
$C_j \in {\calO_{an}}^{r,r}$ are $r\times r$ matrices of holomorphic
functions and one has vanishing commutators
\mbox{$[\del_i-C_i,\del_j-C_j]$}.  Hence the right Koszul complex
$K_\spot(\del-C,D_{an})$ on the operators $\del_i-C_i$ is a free
resolution of $N_{an}$ and the holonomic dual is computed from this
Koszul complex.  In particular, $\DD(N_{an})={D_{an}}^r/\ideal{\del_j
+ C_j \mid j = 1,\ldots,n}$ and so $\rk(N)=\rk(\DD(N))=r$.

The automorphism $x \mapsto -x$ on $\CC^n$ induces, via $\del\mapsto
-\del$, auto-equivalences $N \mapsto N^-$ on the categories of
$R$-modules and of $D$-modules (but not of the category of
$\ZZ^n$-graded \mbox{$\SA$-modules}, since $I_A$ is not preserved
under $(\blank)^-$ unless it is projective).  The formation of
Euler--Koszul complexes is equivariant under this sign change since
$E_i-\beta_i = (E_i-\beta_i)^-$.  Moreover, for $\ZZ^n$-graded
$R$-modules~$N$ we have $D \otimes_{R} N \cong \tau(N\otimes_{R} D)^-$
as $D$-modules, where the tensor products exploit the two different
$R$-structures on $D$.

The ordinary Koszul complex $K_\spot(\boldy;T)$ on a sequence~$\boldy
= y_1,\ldots,y_d$ in a commutative ring~$T$ is isomorphic to its dual
$K^\spot(\boldy;T) = \hom_T(K_\spot(\boldy;T),T)$.  In fact,
$K_\spot(\boldy;T)$ equals $K^\spot(\boldy;T)$ after replacing each
lowered homological index~$i$ by the raised cohomological
index~\mbox{$d-i$} and a suitable sign change in the differentials.
Administering this sign- and index-change to $\calK_\spot(E-\beta;N)$
in Definition \ref{d:ek} yields $\calK^\spot(E-\beta;N)$, whose
cohomology we call the {\em Euler--Koszul cohomology}\/ of~$N$; we
have $\calH^i(E-\beta;N) \cong \calH_{d-i}(E-\beta;N)$.

We shall
apply $\calK_\spot(E-\beta;\blank)$ and $\calK^\spot(E-\beta;\blank)$
to $\ZZ^d$-graded complexes of $R$-modules.  Our conventions for
indexing the resulting double complexes are set up as follows, so that
all homological and cohomological indices are positive.  If
\mbox{$\calF_\spot\colon \calF_0 \from \calF_1 \from \cdots \from
\calF_n$} is a complex of $R$-modules with decreasing lowered indices,
then we write the differentials of the double complex
$\calK_\spot(E-\beta;\calF_\spot)$ pointing downward (the
$\calK_\spot\,$ direction) and to the left (the $\calF_\spot\,$
direction).  On the other hand, if
\mbox{$\calF^\spot\colon \calF^0 \to \calF^1 \to \cdots \to \calF^n$}
has increasing raised indices, then we write the double complex
$\calK^\spot(E-\beta;\calF^\spot)$ with differentials pointing upward
and to the right.
Applying $\hom_T(\blank,T)$ to a complex of free $T$-modules
with decreasing lowered indices (with $T = R$ or $T = D$) yields a
complex with increasing raised indices.

Our theory of toric modules applies to the modules $\ext^i_R(M,R)$
whenever $M$ is toric.

\begin{lemma} \label{l:ext}
If $M$ is a toric $R$-module then  $\ext^i_{R}(M,R)$ is toric
for all~$i$.
\end{lemma}
\begin{proof}
The module $M$ has a composition chain by $\ZZ^d$-graded
$\SA$-modules.  If $M$ is of toric length~$1$ then the finitely
generated $\SA$-module $\ext^i_R(M,R)$ is toric by
Example~\ref{ex:toric}.

In the general case, argue by induction on the toric length of $M$.
The $\ZZ^d$-graded long exact sequence of~$\ext_R(\blank,R)$ arising
from $M_1 \into M \onto M/M_1$ places $\ext^i_R(M,R)$ between two
modules both of which are toric by induction.
\end{proof}

\begin{defn}
Let $\vea$ be the sum $\sum_{j=1}^n a_j$ of the columns of~$A$.  Using
our convention that $\deg(\del_j)=-a_j$, the {\em canonical module}\/
of~$R$ is \mbox{$\omega_R = R(\vea)$}.
\end{defn}

\begin{thm} \label{thm:tot}
If~$M$ is a toric $R$-module, then there is a spectral
sequence
\begin{eqnarray*}
  E_2^{p,q} = \calH^q\bigl(E+\beta;\ext^p_R(M,\omega_R)\bigr)(-\vea)
  &\Longrightarrow&
  H^{p+q}\DD\bigl(\calH_{p+q-n}(E-\beta;M)\bigr){}^-.
\end{eqnarray*}
\end{thm}
\begin{proof}
Let $\calF_\spot$ be a minimal $\ZZ^d$-graded $R$-free resolution
of~$M$.  Consider the double complex
$\tau\!\hom_D(\calK_\spot(E-\beta;\calF_\spot),D)^-$.  With
$\calF^\spot = \hom_R(\calF_\spot,R)$, we claim~that
\begin{eqnarray*}
  \tau\!\hom_D\bigl(\calK_\spot(E-\beta;\calF_\spot),D\bigr){}^-
  &\cong& \calK^\spot(-E-\beta-\vea;\calF^\spot).
\end{eqnarray*}
To see why, begin by noting that $\tau(\blank)^-$ is an isomorphism
when $D$-modules are regarded as $R$-modules, but identifies (for
example) the left $D$-module $D \otimes_R \SA$ with the right
$D$-module $\SA\otimes _R D$. Hence each row of
\mbox{$\tau\!\hom_D(\calK_\spot(E-\beta;\calF_\spot),D)^-$} is a
direct sum of complexes of the form
\mbox{$\tau(\calF^\spot\otimes_R D)^-$}
with cohomology $D\otimes_R \ext^i_R(M,R)$.

On the other hand, note that $\tau(E_i-\beta_i)^- =
-E_i-\beta_i-(\vea)_i$.  Therefore each column of
$\tau\!\hom_D(\calK_\spot(E-\beta;\calF_\spot),D)^-$ is an
Euler--Koszul cocomplex induced by \mbox{$-E-\beta-\vea$} on
$\tau(\calF^i\otimes_R D)^-$. We consider the spectral sequences
associated to this double complex.

Taking first the horizontal and then the vertical cohomology of the
double complex, we obtain $\calH^q(-E-\beta-\vea;\ext^p_R(M,R))\cong
\calH^q\bigl(-E-\beta;\ext^p_R(M,\omega_R)\bigr)(-\vea)$.  We now
determine the abutment by reversing the order of taking horizontal and
vertical cohomology.

The natural projection from the total complex
$\tot\calK_\spot(E-\beta;\calF_\spot)$ to the Euler--Koszul complex
$\calK_\spot(E-\beta;M)$ is a quasi-isomorphism because the rows of
$\calK_\spot(E-\beta;\calF_\spot)$ are resolutions for $M$ positioned
in homological degree~$0$.  Since $\calH_k(E-\beta;M)$ is holonomic
for all~$k$ by Proposition~\ref{p:holonomic}, the complex
$\hom_D(\tot\calK_\spot(E-\beta;\calF_\spot),D)$ has cohomology
\begin{eqnarray*}
  H^{p+q}\hom_D\bigl(\tot\calK_\spot(E-\beta;\calF_\spot),D\bigr) &=&
  \ext^n_D\bigl(\calH_{p+q-n}(E-\beta;M),D\bigr).
\end{eqnarray*}
Applying the standard involution~$\tau$ and the auto-equivalence
$(\blank)^-$ yields
\begin{eqnarray*}
  H^{p+q}\DD\bigl(\tot\calK_\spot(E-\beta;\calF_\spot)\bigr){}^-
  &=& H^{p+q}\DD\bigl(\calH_{p+q-n}(E-\beta;M)\bigr){}^-.
\\[-6.78ex]
\end{eqnarray*}
\end{proof}

\medskip
\begin{remark}
It follows from this spectral sequence that if $M$ is a
Cohen--Macaulay toric module of dimension $d$, 
then $\calH_i(E-\beta;M)=0$ for all
$i>0$ since then $\ext^{n-d+i}_R(M,\omega_R)=0$ for all $i>0$. In the
sequel we link the failure of $M$ to be Cohen--Macaulay to the
appearance of non-vanishing $\calH_i(E-\beta;M)$ for suitable $\beta$
and $i>0$. 
\end{remark}

Let $\frakm = \<\del_1,\dots \del_n\>$ be the $\ZZ^d$-graded maximal
ideal of~$R$.  Given a $\ZZ^d$-graded $R$-module~$N$,
its {\em local cohomology modules}\/
\begin{eqnarray*}
  H^i_\frakm(N) &=& \dlim{t} \ext^i_R(R/\frakm^t,N) 
\end{eqnarray*}
{\em supported at\/~$\frakm$} are $\ZZ^d$-graded.  We refer to
\cite{MS04} for details on the $\ZZ^d$-graded aspects of local
cohomology.  By \cite[Section~3.5]{BH93} there is a natural vector
space isomorphism
\begin{eqnarray*}
  \ext^i_R(N,R)_\alpha &\cong&
  \hom_\CC(H^{d-i}_\frakm(N)_{-\alpha+\vea},\CC)
\end{eqnarray*}
called {\em $\ZZ^d$-graded local duality} (see also \cite{Mil02}).

\begin{defn} \label{d:exdegree}
Fix a toric module~$M$.  A degree $\alpha \in \ZZ^d$ such that
$H^i_\frakm(M)_{-\alpha} \neq 0$ (note the minus sign in the
subscript) for some $i \leq d-1$ is called a {\em true exceptional
degree}\/ of~$M$.  If $\beta \in \CC^d$ lies in the Zariski closure of
the set of true exceptional degrees of~$M$, then $\beta$ is an {\em
exceptional quasi-degree}\/~of~$M$.
\end{defn}

Exactness of the Euler--Koszul complex can be expressed in terms of
\mbox{local cohomology}.

\begin{thm} \label{t:qdeg}
The Euler--Koszul homology \mbox{$\calH_i(E-\beta;M)$} of a toric
module~$M$ over $R$ is nonzero for some $i \geq 1$ if and only if
$\beta \in \CC^d$ is an exceptional quasi-degree of~$M$.  More
precisely, if $k$ equals the smallest homological degree~$i$ for which
$-\beta \in \qdeg(H^i_\frakm(M))$ then $\calH_{d-k}(E-\beta;M)$ is
holonomic of nonzero rank while $\calH_i(E-\beta;M)=0$ for $i > d-k$.
\end{thm}
\begin{proof}
The module~$M$ has dimension at most~$d$, because this is true of
every successive quotient in any toric filtration.  Therefore, using
notation as in Theorem~\ref{thm:tot}, the modules
$\ext^p_R(M,\omega_R))$ can only be nonzero when $p \geq n-d$.

First suppose that~$\beta$ is not an exceptional quasi-degree.  This
means precisely that $\beta$ lies outside of
$\qdeg(\ext^p_R(M,\omega_R))$ for all $p \neq n-d$ by local duality.
Proposition~\ref{p:HH} therefore implies that the only column of the
spectral sequence page~$E_2^{p,q}$ in Proposition~\ref{thm:tot} that
can possibly be nonzero is column~\mbox{$p = n-d$}.  Furthermore, the
highest possible row of a nonzero entry in this column is row~\mbox{$q
= d$}.  Since the cohomology of the abutment is only nonzero when $p+q
\geq n$, it follows that $H^{p+q}\DD(\calH_{p+q-n}(E-\beta;M))^-$ can
only be nonzero when $p + q = n$.  Applying the
auto-equivalence~$(\blank)^-$ and taking holonomic duals, we find that
$\calH_i(E-\beta;M)$ can only be nonzero when~\mbox{$i=0$}.

Now suppose that $\beta$ is an exceptional quasi-degree.  Let $k$ be
the smallest cohomological degree~$i$ such that $\beta \in
-\qdeg(H^i_\frakm(M)) = \qdeg(\ext^{n-i}_R(M,\omega_R))$.  By
Proposition~\ref{p:HH} and Lemma~\ref{l:ext},
\mbox{$\calH^d(E+\beta;\ext^{n-k}_R(M,\omega_R))(-\vea)$} is nonzero.
Moreover, all columns to the right of column~\mbox{$n-k$} and all
rows above row~$d$ in the spectral sequence page $E_2^{p,q}$ of
Theorem~\ref{thm:tot} are zero.  Hence
\begin{eqnarray*}
  \calH^d(E+\beta;\ext^{n-k}_R(M,\omega_R))(-\vea) &=&
  H^{d-k+n}\DD\bigl(\calH_{d-k}(E-\beta;M)\bigr){}^-
\end{eqnarray*}
is by Proposition~\ref{p:holonomic} a holonomic module of nonzero rank.
After applying the auto-equi\-valence~$(\blank)^-$ and taking
holonomic duals, we find that $\calH_{d-k}(E-\beta;M)$ is a holonomic
module of nonzero rank, and the highest nonzero Euler--Koszul homology
of~$M$.%
\end{proof}

\section{Global Euler--Koszul homology as a holonomic family}

\label{sec:global}

In this section we explore the interactions of the Euler--Koszul
functor with our notion of holonomic family.  The main result,
Theorem~\ref{t:fg}, is that the family of $D$-modules
\mbox{$\calH_0(E-\beta;\calM)$} for varying $\beta \in \CC^d$
constitutes a holonomic family over~$\CC^d$.

Resume the notation from Example~\ref{ex:CCd} and
Definition~\ref{d:ek}.  If we give all of the variables $b =
b_1,\ldots,b_d$ degree zero, then the polynomial ring $D[b]$ over the
Weyl algebra~$D$ is $\ZZ^d$-graded, as is its commutative
subalgebra~$R[b]$.

\begin{defn}
Let~$\calM$ be a $\ZZ^d$-graded $D[b]$-module.  The {\em global Euler
endomorphisms}\/ of~$\calM$ are the commuting endomorphisms
$E_1-b_1,\ldots,E_d-b_d$, the $i^\th$ of which acts by
\begin{eqnarray*}
  E_i-b_i:\ m \mapsto (E_i-b_i-\alpha_i)m &\hbox{whenever}& m \in
  \calM_\alpha.
\end{eqnarray*}
Each $\ZZ^d$-graded $R$-module~$N$ yields a {\em global
Euler--Koszul~complex}\/ of left $D[b]$-modules,
\begin{eqnarray*}
  \calK_\spot(E-b;N) &=& \calK_\spot(E-b;D[b] \otimes_{R} N),
\end{eqnarray*}
with $D[b]$-linear homomorphisms.  Write the homology as
$\calH_i(E-b;N) = H_i\calK_\spot(E-b;N)$.
\end{defn}

For any fixed $\beta \in \CC^d$, the zeroth Euler--Koszul homology
$\calH_0(E-\beta;N)$ of~$N$ can be recovered as a fiber (in the sense
of Section~\ref{sec:upper}) of the zeroth global Euler--Koszul
homology $\calH_0(E-b;N)$ of~$N$.  The precise statement, as follows,
is immediate from the definitions.

\begin{lemma} \label{l:fiber}
Suppose that~$N$ is a $\ZZ^d$-graded $R$-module, and that $\beta \in
\CC^d$.  The fiber of $\calM = \calH_0(E-b;N)$ over~$\beta$ is
$\calM_\beta = \calH_0(E-\beta;N)$.\qed
\end{lemma}

In view of Lemma~\ref{l:fiber}, to prove holonomicity for the family
$\calH_0(E-b;M)$ determined by a toric $R$-module~$M$ we need to
establish the coherence condition in Definition~\ref{d:hol}.  For
this, we shall use the criterion in Proposition~\ref{p:filter}.  But
before we can do so, we need to know that Euler operators form
sequences of parameters on the quotients of polynomial rings by
certain initial ideals of toric ideals.  We will get a handle on these
initial ideals via their minimal primes.  For notation, consider the
polynomial ring~$\CC[\xi](x)$ with its usual $\ZZ$-grading, so that
each variable in the list $\xi = \xi_1,\ldots,\xi_n$ has degree~$1$.
By a {\em minimal prime}\/ of a $\CC[\xi](x)$-module~$N$, we mean a
prime minimal among all of those \mbox{containing the annihilator
of~$N$}.

\begin{lemma} \label{l:sysparam}
Let~$N$ be a finitely generated $\ZZ$-graded\/ $\CC[\xi](x)$-module,
and let\/~$\boldy$ be a sequence of $\ZZ$-graded elements
in~$\CC[\xi](x)$.  If\/~\mbox{$\CC[\xi](x)/(\<\boldy\> + \pp)$}
is finite-dimensional as a vector space over\/~$\CC(x)$ for every
minimal prime~$\pp$ of~$N$, then so is~$N/\boldy N$.
\end{lemma}
\begin{proof}
Let $I \subseteq \CC[\xi](x)$ be the annihilator of~$N$.  As~$N$ is
finitely generated over~$\CC[\xi](x)/I$, it is enough to show that
$\<\boldy\> + I$ contains a power of the $\ZZ$-graded maximal ideal
$\frakm = \<\xi\>$.

Suppose that $I$ has radical~$J$.  Then 
if $\frakm^s\subseteq \<\boldy\> + J$  then
$\frakm^{rs}\subseteq (\<\boldy\> +J)^r$, and so 
$\frakm^{rs}\subseteq (\<\boldy\> + J)^r \subseteq \<\boldy\> + J^r
\subseteq \<\boldy\> + I$ for some~$r$.  Since the minimal primes of
$I$ and $J$ are 
identical it is sufficient to consider $J$ instead of $I$.

Now let $J = \pp_1 \cap \cdots \cap \pp_m$ be a primary decomposition
where of course each $\pp_j$ is prime.  Then $\boldy + J = \boldy +
(\pp_1 \cap \pp_2 \cdots \cap \pp_m)$ contains the product
$(\boldy+\pp_1) \cdots (\boldy+\pp_m)$.  Each of the factors
$(\boldy+\pp_j)$ contains a power of~$\frakm$ by the hypothesis, 
so $J+\boldy$ contains the product of these powers of~$\frakm$.%
\end{proof}

We now return to the order filtration on the Weyl algebra~$D$, and the
associated graded ring~$\CC[x,\xi]$ (see Section~\ref{sec:criteria}).
This filtration, when restricted to the subring~$R$ of~$D$, gives an
associated graded ring~$\CC[\xi]$.  It induces on~$R$ the partial
ordering by total degree in which the {\em initial form}\/ of a
polynomial $f \in R$ is the sum~$\IN(f)\in\CC[\xi]$ of all terms of
highest total degree, with the variables $\del_i$ changed to $\xi_i$.
The {\em initial ideal}\/ of any ideal $J \subseteq R$ is the ideal
$\IN(J) = \<\IN(f) \mid f \in J\>$ generated by the initial forms of
all polynomials in~$J$.

In the coming statement, we consider the ideal~$I_F$ from
Notation~\ref{n:toric}.  Let $Ax\xi \subset \CC[\xi](x)$ be the
sequence $\IN(E_1),\ldots,\IN(E_d)$ obtained from $E_1,\ldots,E_d$ by
replacing each~$\del_i$ with~$\xi_i$.

\begin{prop} \label{p:IN}
If\/~\mbox{$\IN(I_F) \subseteq \CC[\xi]$} is the initial ideal
of~$I_F$ in the above sense, then the $\CC(x)$-vector space
$\CC[\xi](x)/\ideal{\IN(I_F), Ax\xi}$ is finite-dimensional over
$\CC(x)$. 
\end{prop}
\begin{proof}
This statement immediately reduces to the case where $I_F = I_A$,
after first replacing~$\CC[\xi](x)$ with $\CC[\xi_j \mid a_j \in
F](x_j \mid a_j \in F)$, and then using any matrix~$A_F'$ composed of
$\dim(F)$ many linearly independent rows of~$F$ to play the role
of~$A$.  We wish to apply Lemma~\ref{l:sysparam}, so we need to
describe the minimal primes of~$\CC[\xi]/\IN(I_A)$.

Define the $(d+1) \times (n+1)$ matrix $\hat A$ by placing a row
$(1,\ldots,1)$ across the top of~$A$, and subsequently adding a
leftmost column $(1,0,0,\ldots,0)$.  If $\xi_0$ is a new variable and
$\hat\xi = \{\xi_0\} \cup \xi$, then $\CC[\xi]/\IN(I_A) \cong
\CC[\hat\xi]/\ideal{I_{\hat A} , \xi_0}$.  Since $\ideal{I_{\hat
A},\xi_0}$ is $\ZZ\hat A$-graded, every minimal prime of
$\CC[\xi]/\IN(I_A)$ is the image in $\CC[\xi] =
\CC[\hat\xi]/\<\xi_0\>$ of a prime ideal~$I_{\hat F} \subset
\CC[\hat\xi]$ for some face~$\hat F$ of~$\NN \hat A$ with $\xi_0 \in
I_{\hat F}$.

Thus we only need that for all faces $\hat F$ of $\NN\hat A$ with
$\xi_0\in I_{\hat F}$, the ring $\CC[\hat\xi](\hat x)/\ideal{I_{\hat
F},Ax\xi}$ has finite dimension as a $\CC(\hat x)$-vector space.  Pick
one such $\hat F$.  Then $\xi_0 \not\in \hat F$, since $\xi_0\in
I_{\hat F}$.  Hence the column indices of~$\hat F$ form a face~$F$
of~$A$ whose $\CC$-linear span does not contain the origin; indeed,
there is a bijection between the faces of~$\hat A$ not containing the
column~$\hat{a_0}$ and the faces of~$A$ not containing the origin in
their span.

We infer that the columns of~$A$ constituting~$F$ lie on a hyperplane
in~$\CC^{d}$ off the origin. Therefore the vector
$(1,\ldots,1)\in\CC^{|F|}$ lies in the rowspan of~$F$.  This implies
that
\begin{eqnarray*}
  \CC[\hat\xi](\hat x)/\ideal{I_{\hat F},Ax\xi}&=&
  \CC[\hat\xi](\hat x)/\ideal{I_{\hat F},Fx\xi}\\&=&
  \CC[\hat\xi](\hat x)/\ideal{I_{\hat F},\hat Fx\xi}\\&=&
  \CC[\hat\xi](\hat x)/\ideal{I_{\hat F},\hat A\hat x\hat \xi},
\end{eqnarray*}
where we write $Fx\xi = \{\IN(E^F_1),\ldots,\IN(E^F_d)\}$, $\hat Fx\xi
= \{\sum_{j\in F}x_j\xi_j, \IN(E^F_1),\ldots,\IN(E^F_d)\}$, and
\mbox{$\hat A\hat x\hat \xi = \{\sum_{j=0}^nx_j\xi_j,
\IN(E_1),\ldots,\IN(E_d)\}$}.

We have thus reduced to the projective situation.  This case of the
proposition follows from the proof of
\cite[Theorem~3.9]{Adolphson-Duke94}, which proceeds by showing
precisely that in the projective situation the Krull dimension of
$\CC[x,\xi]/(\IN(I_A) + \<Ax\xi\>)$ equals~$n$.%
\end{proof}

\begin{thm} \label{t:fg}
For any toric $R$-module~$M$, the sheaf~$\tilde\calM$ on~$\CC^d$ whose
global section module is $\calM = \calH_0(E-b;M)$ constitutes a
holonomic family over~$\CC^d$; in other words, $\calM_\beta =
\calH_0(E-\beta;M)$ is holonomic for all $\beta \in \CC^d$, and
$\calM(x)$ is finitely generated over~$\CC[b](x)$.
\end{thm}
\begin{proof}
That the fibers $\calM_\beta = \calH_0(E-\beta;M)$ are holonomic
modules is a consequence of Lemma~\ref{l:fiber} and
Proposition~\ref{p:holonomic}.

For the coherence condition, suppose first that $M = \SF$ for some
face~$F$.  Using order filtrations as above, the graded
$\CC[b](x)$-module associated to~$\calM(x)$ is a quotient of $N =
\CC[\xi][b](x)/\ideal{\IN(I_F), Ax\xi}$, since
\mbox{$\IN(I_F)+\ideal{Ax\xi}=\IN(I_F)+\IN(E-b)\subseteq
\IN(\ideal{I_F,E-b})$}.  As the generators for $\IN(I_F)$ and~$Ax\xi$
do not involve $b$, Proposition~\ref{p:IN} implies that $N$ is
finitely generated over~$\CC[b](x)$.  Hence $\calM(x)$ is finitely
generated by Proposition~\ref{p:filter}.

Now let $M$ be any toric $R$-module.  For $\alpha \in \ZZ^d$ and with
$b'_i = b_i-\alpha_i$ the isomorphism
$$
\begin{array}{rcccl}
  \calH_0(E-b';\SF)(\alpha) &=& \calH_0(E-b+\alpha;\SF)(\alpha) &=&
  \calH_0(E-b;\SF(\alpha))
\end{array}
$$
together with the independence of $\IN(I_F)$ and $Ax\xi$ of $b$ proves
the theorem if $M$ has toric length~$\ell=1$.  For \mbox{$\ell > 1$},
use the short exact sequence \mbox{$M_1 \into M \onto M/M_1$} from a
toric filtration.  The resulting long exact sequence of global
Euler--Koszul homology $\calH_\spot(E-b;\blank)$ tensored with the
flat module $\CC(x)$ places the module $\calH_0(E-b;M)(x)$ between the
two modules $\calH_0(E-b;M_1)(x)$ and \mbox{$\calH_0(E-b;M/M_1)(x)$},
both of which are finitely generated by induction on~$\ell$.%
\end{proof}

\section{Isomorphism of the two homology theories}

\label{sec:beta-euler}

Let $B$ be a reduced parameter scheme and suppose $\beta\in B_\CC$ is
determined by a regular sequence $\boldy$ in $\Gamma(\calO_B,B)$.  In
Corollary \ref{c:cm} we found that the commutative Koszul complex
$K_\spot(\boldy,\calM(x))$ detects rank-jumps in holonomic families
$\calM$ of $D$-modules at $\beta$.  In the special case of generalized
hypergeometric systems associated to toric modules, where in
particular $B=\CC^d$, we also found in Theorem \ref{t:qdeg} that the
noncommutative Euler--Koszul complex $\calK_\spot(E-\beta;M)$ detects
the quasi-degrees where local cohomology of the toric $R$-module $M$
is non-vanishing in cohomological degree less than~$d$.  In this
section we prove that if $\calM = \calH_0(E-b;M)$ is defined through
the toric module $M$ then these two complexes are isomorphic in the
derived category, and in particular have the same homology.

We begin by noting that the global Euler--Koszul complex can be
interpreted as a collection of commutative Koszul complexes, similarly
to Lemma~\ref{l:honest}.

\begin{lemma} \label{l:honest-global}
Let $N$ be a $\ZZ^d$-graded $R$-module and $\alpha \in \ZZ^d$.  The
$\ZZ^d$-graded degree~$\alpha$ part $\calK_\spot(E-b;N)_\alpha$ of the
global Euler--Koszul complex agrees with the ordinary Koszul complex
$K_\spot(E-b-\alpha; (D[b] \otimes_R N)_\alpha)$ constructed over the
polynomial ring $\Theta[b]\subseteq D[b]$.  That is,
\begin{eqnarray*}
  \calK_\spot(E-\beta;N) &=&
  \bigoplus_{\alpha\in\ZZ^d}K_\spot(E-\beta-\alpha;(D[b]\otimes_R
  N)_\alpha),
\end{eqnarray*}
the right hand side being a direct sum of ordinary Koszul complexes of
$\Theta[b]$-modules.\qed
\end{lemma}

Recall from Theorem~\ref{t:fg} that $\calH_0(E-b;M)$ represents a
holonomic family whenever $M$ is a toric $R$-module. The following
statement links the two homological theories that we have studied:
Euler--Koszul homology and Koszul homology of holonomic families.

\begin{thm} \label{t:euler-beta}
Let~$M$ be a toric $R$-module, and consider the holonomic family on
$\CC^d$ with sections $\calM = \calH_0(E-b;M)$.  For each parameter
vector $\beta \in \CC^d$, the Euler--Koszul homology of~$M$ and the
ordinary Koszul homology of~$\calM$ over\/~$\CC[b]$ are isomorphic:
\begin{eqnarray*}
  \calH_i(E-\beta;M) &\cong& H_i(b-\beta;\calM).
\end{eqnarray*}
\end{thm}
\begin{proof}
Write $\calK^E_\spot = \calK_\spot(E-b;M)$ for the global
Euler--Koszul complex, thought of as a column pointing downward with
its bottom at row~$0$, and let $K^b_\spot = K_\spot(b-\beta,\CC[b])$
be the ordinary Koszul complex on $b-\beta$, thought of as a row
pointing leftward toward its end at column~$0$.  Consider the double
complex \mbox{$K^b_\spot \otimes_{\CC[b]} \calK^E_\spot$}.

Taking horizontal $(K^b_\spot)$-homology first leaves only one column,
namely the leftmost column, and that column is the Euler--Koszul
complex $\calK_\spot(E-\beta;M)$.  It follows that
$H_\spot\tot(K^b_\spot \otimes \calK^E_\spot) \cong
\calH_\spot(E-\beta;M)$.

Now we need to check that $H_\spot\tot(K^b_\spot \otimes
\calK^E_\spot) \cong H_\spot(b-\beta;\calM)$.  For this it is enough
to show that taking the vertical $(\calK^E_\spot)$-homology of the
double complex leaves only one row (i.e., the bottom row), for then
that row is the Koszul complex~$K^b_\spot$ tensored over~$\CC[b]$ with
the global Euler--Koszul homology $\calM = \calH_0(E-b;M)$.

All the homomorphisms in $K^b_\spot \otimes \calK^E_\spot$ are
$\ZZ^d$-graded, so it suffices to check the acyclicity of the vertical
differential separately on each $\ZZ^d$-graded component.  But $(D[b]
\otimes_R M)_\alpha = (D\otimes_R M)_\alpha \otimes_\CC \CC[b]$, so
$\big(b_1 - (E_1-\alpha_1), \ldots, b_d - (E_d-\alpha_d)\big) =
-(E-b-\alpha)$ is a regular sequence in~$\Theta[b]$ on the
$\Theta[b]$-module $(D[b] \otimes_R M)_\alpha$.  Now use
Lemma~\ref{l:honest-global}.
\end{proof}

\section{Combinatorics of hypergeometric ranks}

\label{s:rk-vs-vol}

This section contains our results on our motivating problem: a
geometric description of the rank-jumping locus for generalized
hypergeometric systems (Theorem~\ref{t:MM} and
Corollary~\ref{c:slab}), and a characterization of Cohen--Macaulayness
for semigroup rings through the absence of rank-jumps
(Corollary~\ref{c:GKZ}).

\begin{thm} \label{t:MM}
The rank-jumping parameters $\beta \in \CC^d$ for the hypergeometric
holonomic family $\calH_0(E-b;M)$ of a toric module~$M$ are the
exceptional quasi-degrees~of~$M$.
\end{thm}
\begin{proof}
By Theorem~\ref{t:qdeg}, a parameter~$\beta$ is an exceptional
quasi-degree of~$M$ if and only if the Euler--Koszul homology
$\calH_i(E-\beta;M)$ is of nonzero rank for some~\mbox{$i >
0$}.  This occurs if and only if the homology $H_i(b-\beta;\calM)$ of
the hypergeometric holonomic family $\calM = \calH_0(E-b;M)$ is
holonomic of nonzero rank for some $i > 0$, by
Theorem~\ref{t:euler-beta}.  The nonzero rank condition ensures that
the non-vanishing of $H_i(b-\beta;\calM)$ is equivalent to its
non-vanishing after tensoring with~$\CC(x)$ over~$\CC[x]$.  But $\CC(x)
\otimes_{\CC[x]} H_i(b-\beta;\calM) = H_i(b-\beta;\calM(x))$ by
exactness of localization, and this Koszul homology is nonzero if and
only if $\beta$ is rank-jumping by Corollary~\ref{c:cm}.%
\end{proof}

For a finitely generated $R$-module, the vanishing of~$H^i_\frakm(M)$
for all $i < d$ is equivalent to $M$ having depth~$d$.  This
immediately implies the following.

\begin{cor}\label{c:GKZ}
The holonomic family $\calM_\beta^A= \calH_0(E-\beta;\SA)$ has a
nonempty set of rank-jumping parameters $\beta \in \CC^d$ if and only
if $\SA = R/I_A\,$ fails to be Cohen--Macaulay.\qed
\end{cor}

Not only is the rank-jump locus closed in~$\CC^d$ but it has a nice
geometric structure.

\begin{cor} \label{c:slab}
The set of rank-jumping parameters is a finite union of translates of
linear subspaces $\CC F$ generated by faces~$F$ of\/~$A$.  A translate
of\/~$\CC F$ appears if and only if the prime ideal~$I_F$ lies in the
support of\/~$\ext^{n-i}_{R}(M,R)$ for some $i \le d-1$.
\end{cor}
\begin{proof}
The graded Matlis dual of the local cohomology module~$H^i_\frakm(M)$
is, by local duality, a $\ZZ^d$-graded translate of
$\ext^{n-i}_{R}(M,R)$.  Hence the rank-jumping parameters form a
$\ZZ^d$-graded translate of the negative of the set
$\qdeg(\bigoplus_{i=0}^{d-1}\ext^{n-i}_{R}(M,R))$.%
\end{proof}

\begin{porism}
The set of exceptional parameters for a GKZ hypergeometric system
$\calM^A_\beta = \calH_0(E-\beta;\SA)$ has codimension at least~$2$
in~$\CC^d$.
\end{porism}
\begin{proof}
Set $S = \SA$.  By Corollary \ref{c:slab}, we need
$\ext^{n-i}_R(S,R)_\pp = 0$ for all primes~$\pp$ of dimension $d-1$
and all $i \leq d-1$.  For such primes, $R_\pp$ is regular and local
of dimension $n-d+1$, so $\ext^{n-i}_R(S,R)_\pp =
\ext^{n-i}_{R_\pp}(S_\pp,R_\pp) = 0$ unless $i\geq d-1$.  As
$\ext^{n-d+1}_R(S,R)_\pp$ is Matlis dual over~$R_\pp$ to $H^0_{\pp
R_\pp}(S_\pp)$, which vanishes because $S$ is a domain, we are done.%
\end{proof}

\begin{remark} \label{rk:MM}
Corollary~\ref{c:slab}, and more generally Theorem~\ref{t:semicont},
indicates that the dimensions of the holomorphic solution spaces of
hypergeometric systems~$\calM^A_\beta$ behave quite tamely as
functions of~$\beta$.  However, the methods in this article do not
directly involve a study of the variation of the holomorphic solutions
themselves (cf.\ Remark~\ref{rk:sol}), which could in principle be
much worse.  In fact, work of Saito \cite{Saito-Compositio} implies
that the isomorphism classes of $A$-hypergeometric
systems~$\calM^A_\beta$ do not vary at all algebraically with~$\beta$.
The question therefore remains whether the holomorphic solution space
of a hypergeometric system~$\calM^A_\beta$ varies as a function
of~$\beta$ more like the rank of~$\calM^A_\beta$ or more like the
isomorphism class of~$\calM^A_\beta$.
This problem is particularly important for research toward constructing
explicit solutions to $A$-hypergeometric systems.
\end{remark}

\begin{excise}{%
  \comment{UW: As mentioned to Ezra, I have problems with what
  follows. 1.: If $\beta$ is non-jumping then the natural map
  $\calH_0(E-\beta;S)$ to $\calH_0(E-\beta;\bar S)$ need not be an
  isomorphism. For example, let $S$ be Cohen-Macaulay but not
  saturated, and $\beta$ a quasi-degree of $\bar S/S$.
  2.: There is confusion of notation here. The system
  $\calH_0(E-\beta;\bar S)$ that occurs in 8.6 is not the classical
  GKZ system to $\bar S$ and $\beta$. It has the wrong variables in it
  and therefore the wrong Eulers. So [GKZ87] does not tell us about
  the number of solutions. 
  3.: If we write something like 8.7, we ought to mention [SST], 4.5.2
  and 4.6.1.
  4. It is an important, and as far as I know open, problem to
  determine relations between Euler-Koszul homologies induced by the
  module $M$ that happens to be a module over two semigroup rings
  simultaneously. I would like to have a method that writes down
  solutions for our generalized modules in the way GKZ did it for the
  classical one, and prove that we can write as many as the generic
  rank is equal to.
  5.: I appeal for either rewriting or skipping the remainder of the
  paper.}
  
  Finally, we shift our attention to the {\em non}-rank-jumping
  parameters. In this case, the semigroup ring $\SA$ can be replaced
  by its normalization without changing the holonomic rank of the
  associated hypergeometric system, or indeed, without changing the
  system~at~all.  By a {\em generic parameter}\/ we mean one lying
  outside of a fixed proper Zariski closed subset of~$\CC^d$.
  
  \begin{thm}\label{t:normalization}
  Let $\ol \SA$ be the normalization of~$\SA$.  For generic
  parameters~$\beta$, the hypergeometric systems $\calM^A_\beta =
  \calH_0(E-\beta;\SA)$ and $\calH_0(E-\beta;\ol \SA)$ are isomorphic
  as $D$-modules.
  \end{thm}
  \begin{proof}
  Apply the Euler--Koszul functor to the exact sequence
  $$
    0 \to \SA \to \ol \SA \to \ol \SA/\SA \to 0.
  $$
  Since $\ol \SA/\SA$ is a torsion $\ol \SA$-module, its quasi-degrees
  constitute a Zariski closed proper subset of~$\CC^d$ (in fact, a
  finite union of affine subspaces).  For~$\beta$ outside the
  quasi-degrees of~$\ol \SA/\SA$, all of the Euler--Koszul homology
  $\calH_i(E-\beta;\ol \SA/\SA)$ vanishes by Proposition~\ref{p:HH}.
  It therefore follows from the long exact Euler--Koszul homology
  sequence that for such~$\beta$, the $D$-modules
  $\calH_0(E-\beta;\SA)$ and $\calH_0(E-\beta;\ol \SA)$ are
  isomorphic.%
  \end{proof}
  
  As a consequence, here is a combinatorial formula for the generic
  rank of an $A$-hyper\-geometric system.  Denote by $\vol(A)$ the
  integer obtained by taking the product of~$d!$ times the Euclidean
  volume of the convex hull of the columns of~$A$ and the origin.
  
  \begin{cor}
  For generic $\beta$, the $A$-hypergeometric system~$\calM^A_\beta$
  has $\rk(\calM^A_{\beta})=\vol(A)$.
  \end{cor}
  \begin{proof}
  By Theorem~\ref{t:normalization}, we may assume that $\SA$ is
  normal, and therefore Cohen--Macaulay. In this case, the desired
  rank formula is due to Gelfand, Kapranov and Zelevinsky \cite{GKZ89}
  when $\SA$ is $\ZZ$-graded, and to Adolphson \cite{Adolphson-Duke94}
  in the general case.
  \end{proof}
  
  \begin{remark}
  The question of computing the holonomic rank of a hypergeometric
  system at a rank-jumping parameter is mostly open, although it is
  known that $\rk(\calM^A_\beta) - \vol(A)$ can be arbitrarily large
  \cite{MW}.
  \end{remark}
}\end{excise}%
%



\def\cprime{$'$}
\providecommand{\bysame}{\leavevmode\hbox to3em{\hrulefill}\thinspace}
\providecommand{\MR}{\relax\ifhmode\unskip\space\fi MR }
\providecommand{\MRhref}[2]{%
  \href{http://www.ams.org/mathscinet-getitem?mr=#1}{#2}
}
\providecommand{\href}[2]{#2}


\end{document}